\documentclass{article}

\usepackage{amsmath,amsthm,amssymb}
\usepackage{hyperref}
\usepackage{graphicx}

\newcommand{\unit}{\mathbf{1}}

\newcommand{\diag}{\mathop{\mathrm{diag}}}
\newcommand{\sylow}{\mathcal{K}_p}
\newcommand{\sylown}[1][n]{\mathcal{K}_{p, #1}}
\newcommand{\symm}[1][X]{\mathop{\mathrm{Symm}}(#1)}
\newcommand{\arr}{\longrightarrow}
\newcommand{\mat}{\mathsf{M}}
\newcommand{\bim}{\Phi}
\newcommand{\stencil}{\Xi}
\newcommand{\fieldp}{\mathbb{F}_p}
\newcommand{\autxs}{\mathop{\mathrm{Aut}}(X^*)}
\newcommand{\fautxs}{\mathop{\mathrm{FAut}}(X^*)}

\newtheorem{theorem}{Theorem}[section]
\newtheorem{lemma}[theorem]{Lemma}
\newtheorem{proposition}[theorem]{Proposition}
\newtheorem{corollary}[theorem]{Corollary}

\theoremstyle{definition}

\newtheorem{defi}[theorem]{Definition}
\newtheorem{examp}[theorem]{Example}

\title{Self-similar groups, automatic sequences, and unitriangular representations}
\author{R. Grigorchuk, Y. Leonov, V. Nekrashevych, V. Sushchansky}

\begin{document}
\maketitle

\tableofcontents

\section{Introduction}

Self-similar groups is an active topic of modern group
theory. They initially appeared as interesting examples of groups with
unusual properties
(see~\cite{al:burn_en,susz:burn_en,grigorchuk:80_en,gupta-sidkigroup}). The
main techniques of the theory were developed for the study of these
examples. Later a connection to dynamical systems was discovered
(see~\cite{nek:imgdan_en,nek:book}) via the notion of the iterated
monodromy group. Many interesting problems were solved using
self-similar groups (see~\cite{grigorchuk:notEG,glsz:atiyah,barthvirag,bartnek:rabbit}).

One of the ways to define self-similar groups is to say that they are
groups generated by all states of an automaton (of Mealy type, also
called a transducer, or sequential machine). Especially important case
is when the group is generated by the states of a finite
automaton. All examples mentioned above (including the iterated
monodromy groups of expanding dynamical systems) are like that.

The main goal of this article is to indicate a new relation between
self-similar groups and another classical notion of automaticity:
automatic sequences and matrices. See the 
monographs~\cite{allouchshallit:sequences,haesler:automatic} for
theory of automatic sequences and applications.

More precisely, we are going to study natural linear representations
of self-similar groups over finite fields, and show that matrices
associated with elements of a group generated by a finite automaton
are automatic.

There are several ways to define automatic sequences and matrices. One can use
Moore automata, substitutions (e.g., Morse-Thue substitution leading to the famous Morse-Thue sequence), or Christol's characterization of automatic
sequences in terms of algebraicity of the generating power series
over a suitable finite
field~\cite[Theorem~12.2.5]{allouchshallit:sequences}.
The theory of automatic sequences is rich and is related to many topics in dynamical
systems, ergodic theory, spectral theory of Schr\"odinger operators,
number theory etc., see~\cite{allouchshallit:sequences,haesler:automatic}.

It is well known that linear groups (that is subgroups of groups of
matrices $GL_N(\Bbbk)$, where
$\Bbbk$ is a field) is quite a restrictive class of groups as the
Tits alternative \cite{tits:linear} holds for them.
Moreover, the group of (finite) upper triangular matrices is solvable,
and the group of upper unitriangular
matrices is nilpotent. In contrast, if one uses infinite
triangular matrices over finite field, one
can  get much more groups. In particular, every countable
residually finite $p$-group can be
embedded into the group of upper uni-triangular matrices over the
finite field $\fieldp$.

We will pay special attention to the case when the constructed
representation is a representation by infinite unitriangular
matrices. One of the 
results of our paper is showing that the natural (and optimal in
certain sense) representation by
uni-triangular matrices constructed 
in~\cite{leoneksu:triangle_en,leonov:pgroups} leads to automatic matrices, if the group is generated
by a finite automaton. In particular, the diagonals of these
uni-triangular matrices are automatic sequences. We study them
separately, in particular, by computing their generating series (as
algebraic functions).

The roots of the subject of our paper go back to
L.~Kaloujnin's results on Sylow $p$-subgroups of the
symmetric group~\cite{kalouj:sylowfin,kalouj:pinfty-cras,kaloujnine:pinfty},
theory of groups
generated by finite
automata~\cite{glush:ata,hor:fin_en,al:burn_en,sush:faut_en,grineksu_en}, 
theory of self-similar groups~\cite{bgn,nek:book}, and group actions on rooted
trees~\cite{grigorchuk:branch,handbook:branch,grineksu_en}. 

Note that study of actions on rooted trees (every self-similar group
is, by definition, an automorphism group of the rooted tree of finite
words over a finite alphabet) is
equivalent to the study of residually finite groups by geometric
means, i.e., via representations of them in
groups of automorphisms of rooted trees. The theory of actions
on rooted trees is quite different from the Bass-Serre theory
\cite{serr:trees} of actions on (unrooted)
trees, and uses different methods and tools.
The important case is when a group is residually finite $p$-group ($p$ prime), i.e.,
is approximated by finite $p$-groups. 
The class of residually finite $p$-groups contains groups with many
remarkable properties. For instance,
Golod's $p$-groups that were constructed in~\cite{golod:nilalgebras}
based on Golod-Shafarevich theorem to 
answer a question of Burnside on the existence of a finitely generated
infinite torsion group are residually 
$p$ groups. 
Other important examples are the first self-similar groups mentioned
at the beginning of this introduction.

At the end of the paper we study a notion of uniseriality which plays
an important role in the study of 
actions of groups on finite
$p$-groups~\cite{ddms,leedhamgreenmckay}. Our analysis is based upon
classical results of L.~Kaloujnine on height of automorphisms of
rooted trees~\cite{kalouj:sylowfin,kaloujnine:pinfty}.
Applications of uniseriality to Lie algebras associated with
self-similar groups were, for instance, demonstrated in~\cite{bartholdi-g:lie}.
Proposition~\ref{pr:uniserialitycriterion} gives a simple criterion of
uniseriality of action of a group on 
rooted trees and allows to substitute lemma 5.2  from~\cite{bartholdi-g:lie}. A number of examples is
presented which demonstrate the basic notions, ideas, and results.

\subsection*{Acknowledgments}
The research of the first and third named authors was supported by NSF grants
DMS1207699 and DMS1006280, respectively.

\section{Groups acting on rooted trees}

\subsection{Rooted trees and their automorphisms}
\label{ss:rootedtrees}
Let $X$ be a finite alphabet. Denote by $X^*$ the set of finite words
over $X$, which we will view as the free monoid generated by $X$. It is a disjoint
union of the sets $X^n$ of words of \emph{length $n$}. We denote the
empty word, which is the only element of $X^0$, by $\emptyset$.
We will write elements of $X^n$ either as words $x_1x_2\ldots x_n$, or
as $n$-tuples $(x_1, x_2, \ldots, x_n)$.

We consider $X^*$ as the set of vertices of a rooted tree, defined as the right Cayley graph
of the free monoid. Namely, two vertices of $X^*$ are connected by an
edge if and only if they are of the form $v$, $vx$ for $v\in X^*$ and
$x\in X$. The empty word is the root of the tree. For $v\in X^*$, we
consider $vX^*=\{vw\;:\;w\in X^*\}$ as the set of vertices of a
sub-tree with the root $v$.

Denote by $\autxs$ the group of all automorphisms of the rooted tree
$X^*$. Every element of $\autxs$ preserves the levels $X^n$ of the
tree, and for every $g\in\autxs$ beginning of length $n\le m$ of the word
$g(x_1x_2\ldots x_m)$ is equal to $g(x_1x_2\ldots x_n)$. It follows
that for every $v\in X^*$ the transformation $\alpha_{g, v}:X\arr X$
defined by $g(vx)=g(v)\alpha_{g, v}(x)$ is a permutation, and that the
action of $g$ on $X^*$ is determined by these permutations according
to the rule
\begin{equation}
\label{eq:portaction}
g(x_1x_2\ldots
x_n)=\alpha_{g, \emptyset}(x_1)\alpha_{g, x_1}(x_2)\alpha_{g,
 x_1x_2}(x_3)\ldots\alpha_{g, x_1x_2\ldots x_{n-1}}(x_n).
\end{equation}
The map $v\mapsto\alpha_{g, v}$ from $X^*$ to the symmetric group
$\symm$ is called the \emph{portrait} of the automorphism $g$.

Equivalently, we can represent $g$ by the sequence
\[\tau=[\tau_0, \tau_1, \tau_2, \ldots]\]
of maps $\tau_n:X^n\arr\symm$, where $\tau_n(v)=\alpha_{g, v}$. Such
sequences are called, following L.~Kaloujnine~\cite{kalouj:pinfty-cras},
\emph{tableaux}, and is denoted $\tau(g)$.

If $[\tau_0, \tau_1, \ldots]$ and $[\sigma_0, \sigma_1, \ldots]$ are
tableaux of elements $g_1, g_2\in\autxs$, respectively, then tableau
of their product $g_1g_2$ is the sequence of functions
\begin{equation}
\label{eq:tablproduct}
\tau(g_1g_2)=[\tau_n(g(x_1x_2\ldots x_n))\cdot\sigma_n(x_1x_2\ldots x_n)]_{i=0}^\infty.
\end{equation}
Tableau of the inverse of the element $g_1$ is
\begin{equation}
\label{eq:tablinverse}
\tau(g_1^{-1})=[\tau_n(g^{-1}(x_1x_2\ldots x_n))^{-1}]_{n=0}^\infty.
\end{equation}

Here, and in most of our paper, (except when we will talk about
bisets, i.e., about sets with left and right actions)
group elements and permutations act from the left.

Denote by $X^{[n]}$ the finite sub-tree of $X^*$ spanned by the set of
vertices $\bigcup_{k=0}^nX^k$. The group $\autxs$ acts on $X^{[n]}$,
and the kernel of the action coincides with the kernel of the action
on $X^n$. The quotient of $\autxs$ by the kernel of the action is a
finite group, which is naturally identified with the full automorphism
group of the tree $X^{[n]}$. We will denote this finite group by
$\mathop{\mathrm{Aut}}(X^{[n]})$.

The group $\autxs$ is naturally isomorphic to the inverse limit of the
groups $\mathop{\mathrm{Aut}}(X^{[n]})$ (with respect to the
restriction maps). This shows that $\autxs$ is a profinite group. The
basis of neighborhoods of identity of $\autxs$ is the set of kernels
of its action on the levels $X^n$ of the tree $X^*$.

\subsection{Self-similarity of $\autxs$}

Let $g\in\autxs$, and $v\in X^*$. Then there exists an automorphism of
$X^*$, denoted $g|_v$ such that
\[g(vw)=g(v)g|_v(w)\]
for all $w\in X^*$.

We call $g|_v$ the \emph{section} of $g$ at $v$.
The sections obviously have the properties
\begin{equation}\label{eq:sections}
g|_{v_1v_2}=g|_{v_1}|_{v_2},\qquad
(g_1g_2)|_v=g_1|_{g_2(v)}g_2|_v,
\end{equation}
for all $g, g_1, g_2\in\autxs$ and $v, v_1, v_2\in X^*$.

The portrait of the section $g|_v$ is obtained by restricting the
portrait of $g$ onto the subtree $vX^*$, and then identifying $vX^*$
with $X^*$ by the map $vw\mapsto w$.

\begin{defi}
The set $g|_X=\{g|_v\;:\;v\in X^*\}\subset\autxs$ for $g\in\autxs$, is
called the \emph{set of states} of $g$. An automorphism $g\in\autxs$
is said to be \emph{finite state} if $g|_X$ is finite.
\end{defi}

It follows from~\eqref{eq:sections} that
\[g^{-1}|_X=(g|_X)^{-1},\qquad (g_1g_2)|_X\subset g_1|_Xg_2|_X,\]
which implies that the set of finite state elements of $\autxs$ is a
group. We call it the \emph{group of finite automata}, and denote it
$\fautxs$. This name comes
from interpretation of elements of $\autxs$ with automata
(transducers), see~\ref{ss:MMautomata} below. Namely, the set of state of the automaton corresponding
to $g$ is $g|_X$. The element $g\in g|_X$ is the \emph{initial
 state}. If the current state of the automaton is $h$, and it reads
a letter $x\in X$ on its input, then it outputs $h(x)$ and changes it
current state to $h|_x$. It is easy to check that if we give the
consecutive letters of a word $x_1x_2\ldots x_n$ on input of the
automaton with the initial state $g$, then we will get on output the
word $g(x_1x_2\ldots x_n)=g(x_1)g|_{x_1}(x_2)g|_{x_1x_2}(x_3)\ldots$,
compare with~\eqref{eq:portaction}.

Every element $g\in\autxs$ is uniquely determined by the permutation
$\pi$ it defines on the first level $X\subset X^*$ and the first level
sections $g|_x$, $x\in X$. In fact, the map
\[g\mapsto \pi\cdot(g|_x)_{x\in X}\]
is an isomorphism of $\autxs$ with the wreath product
$\symm\ltimes\autxs^X=\symm\wr\autxs$. We call the isomorphism
\[\phi:\autxs\arr\symm\wr\autxs:g\mapsto\pi\cdot(g|_x)_{x\in X}\]
the \emph{wreath recursion}.

For a fixed ordering $x_1, x_2, \ldots, x_d$ of the letters of $X$,
the elements of $\symm\wr\autxs$ are written as $\pi(g_1, g_2, \ldots,
g_d)$, where $\pi\in\symm[d]$ and $g_i=g|_{x_i}$.

\begin{defi}
\label{def:selfsimilar}
A subgroup $G\le\autxs$ is said to be
\emph{self-similar} if $g|_x\in G$ for all $g\in G$ and $x\in X$.
\end{defi}

In other words, a group $G\le\autxs$ is self-similar if restriction of
the wreath recursion onto $G$ is a homomorphism $\phi:G\arr\symm\wr
G$. Note that wreath recursion is usually not an isomorphism (but is
an embedding, since we assume that $G$ acts faithfully on $X^*$).

\begin{examp}
\label{ex:admach}
Let $X=\{0, 1\}$. Consider the automorphism of the tree $X^*$ given by
the rules
\[a(0w)=1w,\qquad a(1w)=0a(w).\]
These rules can be written using wreath recursion as
$\phi(a)=\sigma(1, a)$, where $\sigma=(01)$ is the transposition. We
will usually omit $\phi$, and write
\[a=\sigma(1, a),\]
thus identifying $\autxs$ with $\symm\wr\autxs$.

The automorphism $a$ is called the (binary) \emph{adding machine},
since it describes the process of adding one to a dyadic integer:
$a(x_1x_2\ldots x_n)=y_1y_2\ldots y_n$ if and only if
\[(x_1+2x_2+2^2x_3+\cdots+2^{n-1}x_n)+1=y_1+2y_2+\cdots+2^{n-1}y_n\pmod{2^n}.\]

The group generated by $a$ (which is infinite cyclic)
is self-similar, and is a subgroup of the group of finite automata.
\end{examp}

\begin{examp}
\label{ex:grigorchuk}
Consider the group $G$ generated by the elements $a, b, c, d$ that are
defined inductively by the recursions
\[a=\sigma,\quad b=(a, c),\quad c=(a, d),\quad d=(1, b).\]
Here $\sigma$, as before, is the transposition $(01)$, and when we
omit either the element of $\symm$ or the element of $\autxs^X$ when
writing elements of $\symm\wr\autxs$, we assume that it is equal to
the identity element of the respective group.

The group $G$ is then a self-similar subgroup of the group of finite
automata. It is the Grigorchuk group, defined in~\cite{grigorchuk:80_en}.
\end{examp}

\subsection{Self-similarity bimodule}
\label{ss:ssbimodule}

We can identify the letters $x\in X$ with transformations $v\mapsto
xv$ of the set $X^*$. Then the identity $g(xv)=yh(v)$ for $x\in X$,
$y=g(x)$, and $h=g|_x$
is written as equality of compositions of transformations:
\[g\cdot x=y\cdot h.\]

Consider the set $X\cdot G$ of compositions of the form $x\cdot g$,
i.e., transformations $v\mapsto xg(v)$, $v\in X^*$. It is closed with respect
to pre- and post-compositions with the elements of $G$:
\[(x\cdot g)\cdot h=x\cdot (gh),\qquad h\cdot (x\cdot g)=h(x)\cdot
(h|_xg).\]
We get in this way a \emph{biset}, i.e., a set with two commuting left and
right actions of the group $G$.

Let $\Bbbk$ be a field, and let $\Bbbk[G]$ be the group ring over
$\Bbbk$. Denote by $\bim$ the vector space $\Bbbk^{X\cdot G}$ spanned
by $X\cdot G$. Then the left and the right actions of $G$ on $X\cdot
G$ are extended by linearity to a structure of a $\Bbbk[G]$-bimodule
on $\bim$. We will denote by ${}_G\bim$ and $\bim_G$
the space $\bim$ seen as a left and a right $\Bbbk[G]$-module, respectively.

It follows directly from the definition of the right action of $G$ on
$X\cdot G$ that $X$ (identified with $X\cdot 1$) is a free basis of
$\bim_G$. The left action is not free in general, since it is possible
to have $g(xv)=xv$ for all $v\in X^*$ and for a non-trivial element
$g\in G$, which will imply, by definition of the left action, that
$g\cdot x=x$.

For every element $a\in\Bbbk[G]$ the map $v\mapsto a\cdot v$
for $v\in\bim$ is an endomorphism of $\bim_G$, denoted
$\stencil(a)$. The map
$\stencil:\Bbbk[G]\arr\mathop{\mathrm{End}}\bim_G$ is
obviously a homomorphism of $\Bbbk$-algebras.

After fixing a basis of the right module $\bim_G$ (for example
$X$), we can identify the algebra of endomorphisms
$\mathop{\mathrm{End}}\bim_G$ of the right $\Bbbk[G]$-module
$\bim_G$ with the algebra of $|X|\times|X|$ matrices over
$\Bbbk[G]$. In this case the homomorphism $\stencil:\Bbbk[G]\arr
\mat_{|X|}(\Bbbk[G])\cong\mathop{\mathrm{End}}\bim_G$
is called the \emph{matrix recursion} associated with the self-similar
group $G$ (and the basis of the right module).

More explicitly, if $B=\{e_1, \ldots, e_d\}$ is a basis of the right $\Bbbk[G]$-module
$\bim_G$, then, for $a\in\Bbbk[G]$ the matrix $\stencil(a)=(a_{i,
 j})_{1\le i, j\le d}\in\mat_d(\Bbbk[G])$ is given by the condition
\[a\cdot e_j=\sum_{i=1}^de_i\cdot a_{i, j}.\]

If we use the basis $\{x_1, x_2, \ldots, x_d\}=X$ of the right module
$\bim_G$, then the matrix recursion $\stencil$ is a direct rewriting
of the wreath recursion $\phi:G\arr\symm\wr G$ in matrix
terms. Namely, $\stencil(g)$ is the matrix with entries $a_{ij}$,
$1\le i, j\le d$, given by the rule
\begin{equation}\label{eq:matrecursion}
a_{ij}=\left\{\begin{array}{rl}g|_{x_j} & \text{if $g(x_j)=x_i$,}\\
  0 & \text{otherwise}\end{array}\right.
\end{equation}

\begin{examp}
\label{ex:admach2}
The adding machine recursion $a=\sigma(1, a)$ is defined in the terms
of the bimodules as
\[a\cdot x_0=x_1,\qquad a\cdot x_1=x_0\cdot a,\]
where $x_0, x_1$ are identified with the symbols $0, 1$, respectively,
from Example~\ref{ex:admach}.

It follows that the recursion is written in matrix
form as
\[\stencil(a)=\left(\begin{array}{cc} 0 & a\\ 1 &
  0\end{array}\right).\]

The recursive definition of the generators $a, b, c, d$ of the
Grigorchuk group is written as
\[\stencil(a)=\left(\begin{array}{cc} 0 & 1\\ 1 & 0\end{array}\right),
\stencil(b)=\left(\begin{array}{cc} a & 0\\ 0 & c\end{array}\right),
\stencil(c)=\left(\begin{array}{cc} a & 0\\ 0 & d\end{array}\right),
\stencil(d)=\left(\begin{array}{cc} 1 & 0\\ 0 & b\end{array}\right).\]
\end{examp}

When we change the basis of the right module $\bim_G$, we
just conjugate the map $\stencil$ by the transition matrix. Namely, if
$\{x_1, \ldots, x_d\}$ and $\{y_1, \ldots, y_d\}$ are bases of the
right module $\bim_G$, then we can write $y_j=\sum_{i=1}^\infty
y_i\cdot b_{i, j}$ for $b_{i, j}\in\Bbbk[G]$. Then the matrix $T=(b_{i,
 j})_{1\le i, j\le d}$ is the transition matrix from the basis
$\{x_i\}_{1\le i\le d}$ to the basis $\{y_i\}_{1\le i\le d}$.

\begin{examp}
\label{ex:admach3}
Consider again the adding machine example. Let us take, instead of the
standard basis $\{x_0, x_1\}=X$, the basis $y_0=x_0+x_1$,
$y_1=x_1$. (Here we replace the letters $0, 1$ of the binary alphabet
by $x_0$ and $x_1$, respectively, in order not to confuse them with
elements $0, 1\in\Bbbk[G]$.) Then the transition matrix to the new basis is
$T=\left(\begin{array}{cc} 1 & 0\\ 1 & 1\end{array}\right)$. It
inverse is $T^{-1}=\left(\begin{array}{cc}1 & 0\\ -1 &
  1\end{array}\right)$. Consequently, the matrix recursion in the
new basis is
\[a\mapsto T^{-1}\left(\begin{array}{cc} 0 & a\\ 1 &
  0\end{array}\right)T=\left(\begin{array}{cc} a & a\\ 1-a &
  -a\end{array}\right).\]
This can be checked directly:
\[a\cdot y_0=a\cdot (x_0+x_1)=x_1+x_0\cdot
a=y_1+(y_0-y_1)\cdot a=y_0\cdot a+y_1\cdot(1-a)\]
and
\[a\cdot y_1=a\cdot x_1=x_0\cdot a=y_0\cdot a-y_1\cdot a.\]
If we take the basis $\{y_0=x_0, y_1=x_1\cdot a\}$, then matrix
recursion becomes
\[a\mapsto\left(\begin{array}{cc} 1 & 0\\ 0 &
  a^{-1}\end{array}\right)\cdot
\left(\begin{array}{cc} 0 & a\\ 1 & 0\end{array}\right)\cdot
\left(\begin{array}{cc} 1 & 0 \\ 0 & a\end{array}\right)=
\left(\begin{array}{cc} 0 & a^2\\ a^{-1} & 0\end{array}\right).\]
If the basis is a subset of $X\cdot G$, then the matrix
recursion corresponds to a wreath recursion $G\mapsto \symm\wr G$. For
instance, in the last example the matrix recursion corresponds to the wreath recursion
\[a\mapsto\sigma(a^{-1}, a^2).\]
This wreath recursion describes the process of adding $1$ to a dyadic
numbers in the binary numeration system with digits $0$ and $3$. For
more on changes of bases in the biset $X\cdot G$ and the corresponding
transformations of the wreath recursion see~\cite{nek:book,nek:filling}.
\end{examp}

If $\bim_1$ and $\bim_2$ are bimodules over a $\Bbbk$-algebra $A$,
then their tensor product $\bim_1\otimes\bim_2$ is the quotient of the
$\Bbbk$-vector space spanned by $\bim_1\times\bim_2$ by the sub-space
generated by the elements of the form
\[(v_1, a\cdot v_2)-(v_1\cdot a, v_2)\]
for $v_1\in\bim_1$, $v_2\in\bim_2$, $a\in A$. It is a
$\Bbbk[G]$-bimodule with respect to the actions $a\cdot (v_1\otimes
v_2)=(a\cdot v_1)\otimes v_2$ and $(v_1\otimes v_2)\cdot a=v_1\otimes
(v_2\cdot a)$.

If $\bim_2$ is a left
$A$-module, and $\bim_1$ is an $A$-bimodule, then the left module
$\bim_1\otimes\bim_2$ is defined in the same way.

Let $\bim$, as above, be the bimodule associated with a self-similar
group $G$. Then $X$ is a basis of the right $\Bbbk[G]$-module $\bim_G$, and the set
\[X^n=\{x_1\otimes\cdots\otimes x_n\;:\;x_i\in X\}\]
is a basis of the right module $\bim^{\otimes n}_G$, which is hence a
free module. Note that
$X^n\cdot G$ is the basis of $\bim^{\otimes n}$ as a vector space over
$\Bbbk$.

We identify $x_1\otimes\cdots\otimes x_n$ with the word $x_1\cdots x_n$.
The left module structure on $\bim^{\otimes n}$ is given by the rules similar to the definition
of $\bim$:
\begin{equation}\label{eq:nbim}
g\cdot v=g(v)\cdot g|_v
\end{equation}
for $v\in X^n$ and $g\in G$. In particular, up to an ordering of the
basis $X^n$, the associated matrix recursion $\stencil^n:\Bbbk[G]\arr
\mat_{|X^n|}(\Bbbk[G])$ is obtained from the recursion
$\stencil^{n-1}:\Bbbk[G]\arr\mat_{|X^{n-1}|}(\Bbbk[G])$
by replacing every entry $a_{ij}$ of the matrix $\stencil^{n-1}(a)$ by
the matrix $\stencil(a_{ij})$.

\begin{examp}
The matrix recursion $G\arr\mat_4(\Bbbk[G])$ for the adding machine
(in the standard basis $X^2$) is
\[a\mapsto\left(\begin{array}{cc|cc}0 & 0 & 0 & a\\ 0 & 0 & 1 & 0\\
  \hline 1 & 0 & 0 & 0\\ 0 & 1 & 0 & 0\end{array}\right),\]
which is obtained by iterating the matrix recursion
\[a\mapsto\left(\begin{array}{cc}0 & a\\ 1 & 0\end{array}\right).\]
In this case the basis $X^2$ is ordered in the lexicographic order
$x_0x_0<x_0x_1<x_1x_0<x_1x_1$. But since $a$ is the adding machine,
and it describes adding 1 to a dyadic integer that is written is such a
way that the less significant digits come before the more significant
ones, it is more natural to order the basis in the \emph{inverse}
lexicographic order $x_0x_0<x_1x_0<x_0x_1<x_1x_1$. In this case the
matrix recursion becomes
\[a\mapsto\left(\begin{array}{cccc}0 & 0 & 0 & a\\ 1 & 0 & 0 & 0\\ 0 &
  1 & 0 & 0\\ 0 & 0 & 1 & 0\end{array}\right).\]\end{examp}

\begin{proposition}
\label{pr:Tn}
Let $T$ be the transition matrix from the basis $X$ of
$\bim_G$ to a basis $Y$. Suppose that all entries of $T$ are
elements of $\Bbbk$. Then the transition matrices $T_n$ from
the basis $X^{\otimes n}$ to $Y^{\otimes n}$ is equal to
\[T_n=\underbrace{T\otimes T\otimes\cdots\otimes T}_{|X|},\]
where $\otimes$ is the Kronecker product of matrices.
\end{proposition}

\begin{proof}
Let $T_{n-1}=(a_{u, v})_{u\in X^{\otimes (n-1)}, v\in Y^{\otimes
  (n-1)}}$, i.e.,
\[v=\sum_{u\in X^{\otimes (n-1)}}u\cdot a_{u, v}\] for all $v\in Y^{\otimes
 (n-1)}$. Similarly, denote $T=(a_{x, y})_{x\in X, y\in Y}$.
Then
\begin{multline*}y\otimes v=y\otimes \sum_{u\in X^{\otimes
   (n-1)}}u\cdot a_{u, v}=\left(\sum_{x\in
  X}x\cdot a_{x, y}\right)\otimes\sum_{u\in
 X^{\otimes (n-1)}}u\cdot a_{u, v}=\\ \sum_{u\otimes x\in
 X^{\otimes n}}x\cdot\otimes a_{x, y}\cdot u\cdot a_{u, v}=
\sum_{u\otimes x\in
 X^{\otimes n}}x\otimes u\cdot a_{x, y}\cdot a_{u, v},
\end{multline*}
which shows that
\[a_{xu, yv}=a_{x, y}a_{u, v},\]
which agrees with the definition of the Kronecker product.
\end{proof}

In other words, we can write
\begin{equation}\label{eq:Tn}
T_n=T^{(n)}\left(\begin{array}{cccc}T_{n-1} & 0 & \ldots & 0\\ 0 & T_{n-1} &
  \ldots & 0\\ \vdots & \vdots & \ddots & \vdots\\ 0 & 0 & \ldots & T_{n-1}\end{array}\right),
\end{equation}
where $T_1=T$, and $T^{(n)}$ is the matrix $T$ in which
each entry $a_{ij}$ is replaced by $a_{ij}$ times the unit
matrix of dimension $|X|^{n-1}\times |X|^{n-1}$. Here the
rows and columns of $T_n$ correspond to the elements of $X^{\otimes n}$ and
$Y^{\otimes n}$, respectively, ordered lexicographically.

It is easy to see from the proof, that
in the general case (when not all entries of $T$ are elements of
$\Bbbk$), the formula~\eqref{eq:Tn} remains to be true, if we replace
$T^{(n)}$ by the image of $T$ under the $(n-1)$st iteration of the matrix
recursion (in the basis $X$).

\begin{examp}
\label{ex:Hadamard}
Let $\Bbbk=\mathbb{R}$, $X=\{x_0, x_1\}$. Consider a new basis of $\bim_G$
\[\left\{y_0=\frac{x_0+x_1}{\sqrt{2}},
 y_1=\frac{x_0-x_1}{\sqrt{2}}\right\}.\] The transition matrix to the
new basis is \[\left(\begin{array}{rr} 1/\sqrt{2} & 1/\sqrt{2}\\
  1/\sqrt{2} &
  -1/\sqrt{2}\end{array}\right)=\frac{1}{\sqrt{2}}\left(\begin{array}{rr}
  1 & 1 \\ 1 & -1\end{array}\right).\]

Then the transition matrix from
 $X^{\otimes n}$ to $Y^{\otimes n}$ satisfies the recursion
\[H_n=\frac{1}{\sqrt{2}}\left(\begin{array}{rr}
  1 & 1 \\ 1 & -1\end{array}\right)\left(\begin{array}{cc} H_{n-1} &
  0 \\ 0 & H_{n-1}\end{array}\right)=\frac{1}{\sqrt{2}}\left(\begin{array}{rr}
  H_{n-1} & H_{n-1} \\ H_{n-1} & -H_{n-1}\end{array}\right).\]
\end{examp}

\subsection{Inductive limit of $\Bbbk^{X^n}$}
\label{ss:indlimit}

Let $\Bbbk^{X^n}$ be the vector space of functions $X^n\arr\Bbbk$.
It is naturally isomorphic to the $n$th tensor
power of $\Bbbk^X$. The isomorphism maps an elementary tensor
$f_1\otimes f_2\otimes\cdots\otimes f_n$ to the function
\[f_1\otimes f_2\otimes\cdots\otimes f_n(x_1x_2\ldots x_n)=
f_1(x_1)f_2(x_2)\cdots f_n(x_n).\]
More generally, we have natural isomorphisms
$\Bbbk^{X^n}\otimes\Bbbk^{X^m}\cong\Bbbk^{X^{n+m}}$
defined by the equality
\[f_1\otimes f_2(x_1x_2\ldots x_{n+m})=f_1(x_1x_2\ldots
x_n)f_2(x_{n+1}x_{n+2}\ldots x_{n+m}).\]

We denote by $\delta_v$, for $v\in X^n$ the delta-function of $v$,
i.e., the characteristic function of $\{v\}$. It is an element of
$\Bbbk^{X^n}$. Note that \[\delta_{x_1x_2\ldots
 x_n}=\delta_{x_1}\otimes\delta_{x_2}\otimes\cdots\otimes\delta_{x_n},\]
with respect to the above identification of $\Bbbk^{X^{n+m}}$ with $\Bbbk^{X^n}\otimes\Bbbk^{X^m}$.

Let $G\le\autxs$.
Denote by $\pi_n$ the natural permutational representation of $G$ on $\Bbbk^{X^n}$
coming from the action $G$ on $X^n$. It is given by the rule
$\pi_n(\delta_v)=\delta_{g(v)}$, i.e., by
\[\pi_n(g)(f)(v)=f(g^{-1}(v)),\qquad f\in\Bbbk^X, v\in X^n.\]

Denote by $V_n$ the vector space
$\Bbbk^{X^n}$ seen as a left $\Bbbk[G]$-module of the representation
$\pi_n$, and by $[\varepsilon]$ the left $\Bbbk[G]$-module of the trivial
representation of $G$. More explicitly, it is a one-dimensional vector
space over $\Bbbk$ spanned by an element $\varepsilon$, together with
the left action of $\Bbbk[G]$ given by the rule
\[g\cdot\varepsilon=\varepsilon\]
for all $g\in G$.
The following proposition is a direct corollary of~\eqref{eq:nbim}.

\begin{proposition}
\label{pr:bimeps}
The left module $V_n$ is isomorphic to
$\bim^{\otimes n}\otimes[\varepsilon]$. The isomorphism is
the $\Bbbk$-linear extension of the map
$\delta_{x_1x_2\ldots x_n}\mapsto x_1\otimes x_2\otimes\cdots\otimes
x_n\otimes\varepsilon$ for $x_i\in X$.
\end{proposition}

Denote by $\unit$ the function $\sum_{x\in X}\delta_x\in V_1$ taking
constant value $1\in\Bbbk$. We have then, for every $f\in V_n=\Bbbk^{X^n}$,
\[f\otimes\unit(x_1x_2\ldots x_{n+1})=f(x_1x_2\ldots x_n).\]
The following proposition is straightforward.

\begin{proposition}
The map $\iota_n:v\mapsto v\otimes\unit:V_n\arr V_{n+1}$ is an embedding
of the left
$\Bbbk[G]$-modules. In other words,
\[\pi_{n+1}(g)(f\otimes\unit)=\pi_n(g)(f)\otimes\unit\]
for all $g\in G$ and $f\in V_n$.
\end{proposition}

The space $X^\omega=\{x_1x_2\ldots\;:\;x_i\in X\}$
has a natural topology of a direct (Tikhonoff)
power of the discrete space $X$. A basis of this topology consists of
the cylindrical sets $vX^\omega$, for $v\in X^*$.

Denote by $C(X^\omega, \Bbbk)$ the vector space of maps
$f:X^\omega\arr\Bbbk$ such that $f^{-1}(a)$ is open and closed (\emph{clopen}) for every
$a\in\Bbbk$. In other words, $C(X^\omega, \Bbbk)$ is the space of
all continuous maps $f:X^\omega\arr\Bbbk$, where $\Bbbk$ is taken with
discrete topology. Note that the set of values of any element of
$C(X^\omega, \Bbbk)$ is finite, since $X^\omega$ is compact.

For example, a map
$f:X^\omega\arr\mathbb{R}$ belongs to $C(X^\omega, \mathbb{R})$ if and
only it is continuous and has a finite set of values.

The group $G$ acts naturally on $X^\omega$ by homeomorphisms,
hence it also acts naturally on the space
$C(X^\omega, \Bbbk)$ by the rule
\[g(\xi)(w)=\xi(g^{-1}(w))\]
for $g\in G$, $\xi\in C(X^\omega, \Bbbk)$, and $w\in X^\omega$.

For every $f\in V_n=\Bbbk^{X^n}$ consider the natural extension of
$f:X^n\arr\Bbbk$ to a function on $X^\omega$:
\[f(x_1x_2\ldots)=f(x_1x_2\ldots x_n).\]
For example, the delta-function $\delta_v$ is extended to the
characteristic functions of the subset $vX^\omega$, which we will also
denote $\delta_v$.

It is easy to see that this defines an embedding of
$\Bbbk[G]$-modules $V_n\arr C(X^\omega, \Bbbk)$. Moreover, these
embeddings agree with the embeddings $\iota_n:V_n\arr V_{n+1}$.

Denote by $V_\infty$ the direct limit of the $G$-modules $V_n$ with respect to
the maps $\iota_n$. We will denote by $\pi_\infty$ the corresponding
representation of $G$ on $V_\infty$.

\begin{proposition}
The module $V_\infty$ is naturally isomorphic to the left
$\Bbbk[G]$-module $C(X^\omega, \Bbbk)$.
\end{proposition}

\begin{proof}
The set $\{f^{-1}(t)\;:\;t\in\Bbbk, f^{-1}(t)\ne\emptyset\}$ is a finite covering of
$X^\omega$ by clopen disjoint sets. Every clopen set of $X^\omega$ is a
finite union of cylindrical sets of the form $vX^\omega$, for $v\in
X^*$. Consequently, there exists $n$ such that $f$ is constant on
every cylindrical set of the form $vX^\omega$ for $v\in X^n$. Then
$f\in\Bbbk^{X^n}$ in the identification of $\Bbbk^{X^n}$ with a
subspace of $C(X^\omega, \Bbbk)$, described above. It follows that
the inductive limit of $\Bbbk^{X^n}$ coincides with $C(X^\omega,
\Bbbk)$. We have already seen that the representations $\pi_n$ agree
with the representation of $G$ on $C(X^\omega, \Bbbk)$, restricted
to $V_n=\Bbbk^{X^n}$) which finishes the proof.
\end{proof}

Let $\mathsf{B}$ be a basis of the $\Bbbk$-vector space $\Bbbk^X$ such
that the constant one function $\unit$ belongs to $\mathsf{B}$. Then
$\mathsf{B}^{\otimes n}$ is a basis of the $\Bbbk$-vector space
$\Bbbk^{X^n}=V_n$, and we have $\iota_n(\mathsf{B}^{\otimes n})\subset\mathsf{B}^{\otimes {n+1}}$.
Then the inductive limit $\mathsf{B}_\infty$ of the bases $\mathsf{B}^{\otimes n}$
with respect to the maps $\iota_n$ is a basis of $C(X^\omega, \Bbbk)=V_\infty$. The
elements of this basis are equal to functions of the form
\[f(x_1x_2\ldots)\mapsto f_1(x_1)f_2(x_2)\cdots,\]
where $f_i\in\mathsf{B}$ and all but a finite number of the functions $f_i$ are
equal to the constant one.

\begin{examp}
Suppose that the field $\Bbbk\cong\mathbb{F}_q$ is finite, and let $X=\Bbbk$.
Then the functions $e_k:X\arr\Bbbk:x\mapsto x^k$ for $k=1, 2,\ldots, q-1$ together
with the constant one function $\unit$, formally denoted $x^0$, form a basis of $V_1$.

The corresponding basis of $C(X^\omega, \Bbbk)$ is equal to the set of
all finite monomial functions
\[f(x_1, x_2, \ldots)=x_1^{k_1}x_2^{k_2}\cdots,\]
where all but a finite number of powers $k_i$ are equal to zero.

Writing the elements of $C(X^\omega, \Bbbk)$ in this basis amounts to
representing them as polynomials.
\end{examp}

\begin{examp}
Let $X=\{x_0, x_1\}$, $\mathop{\mathrm{char}}\Bbbk\ne 2$, and let $\mathsf{W}$ be the
basis of $\Bbbk^{X}$ consisting of functions
$y_0=\delta_{x_0}+\delta_{x_1}=\unit$ and
$y_1=\delta_{x_0}-\delta_{x_1}$. The corresponding basis $\mathsf{W}_\infty$ of
$C(X^\omega, \Bbbk)$ is called the \emph{Walsh basis}, see~\cite{walsh}.

For $\Bbbk=\mathbb{C}$, the Walsh basis is an orthonormal set of
complex-valued functions on $X^\omega$ with
respect to the uniform Bernoulli measure on $X^\omega$. This is a
direct corollary of the fact that $\{y_0, y_1\}$ is orthonormal. Since
$\mathsf{W}_\infty$ is a basis of the linear space of continuous functions
$X^\omega\arr\mathbb{C}$ with finite sets of values, and this space is dense in the
Hilbert space $L^2(X^\omega)$, the Walsh basis is an orthonormal basis
of $L^2(X^\omega)$.

We can use Proposition~\ref{pr:Tn} to find transition matrices from
$\{\delta_v\}_{v\in X^n}$ to the
basis $\mathsf{W}^{\otimes n}$ (just use the proposition for the case of the trivial
group $G$). In the case of Walsh basis we get the matrices from Example~\ref{ex:Hadamard},
but without $1/\sqrt{2}$:
\[H_n=\left(\begin{array}{cc} H_{n-1} & H_{n-1} \\ H_{n-1} &
  -H_{n-1}\end{array}\right),\]
compare with Example~\ref{ex:Hadamard}. These matrices are examples of
Hadamard matrices (i.e., matrices whose entries are +1 and -1 and
whose rows are orthogonal) and were constructed for the first time by
J.J.Sylvester~\cite{sylvester:matrices}. They are also called
\emph{Walsh matrices}.

\begin{figure}
\centering
\includegraphics{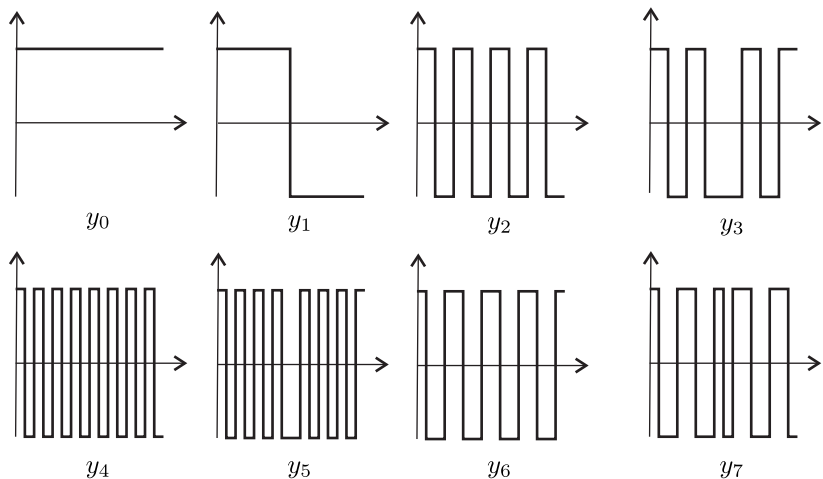}
\caption{Walsh basis}
\label{fig:walsh}
\end{figure}

See Figure~\ref{fig:walsh}, where graphs of
the first eight elements of the Walsh basis are shown. Here we identify $\{0,
1\}^\omega$ with the unit interval $[0, 1]$ via real binary numeration system.
\end{examp}

\begin{examp}
A related basis of $C(X^\omega, \Bbbk)$ is the \emph{Haar basis},
which is constructed in the following way. Again, we assume that
characteristic of $\Bbbk$ is different from 2, and $X=\{x_0,
x_1\}$. Let $y_0=\unit$ and $y_1=\delta_{x_0}-\delta_{x_1}$, as in the
previous example. Let us construct an increasing sequence of bases $\mathsf{Y}_n$
of $\Bbbk^{X^n}<C(X^\omega, \Bbbk)$ in the following way. Let
$\mathsf{Y}_0=\{y_0\}$. Define then inductively:
\[\mathsf{Y}_{n+1}=\mathsf{Y}_n\cup \{\delta_v\otimes y_1\;:\;v\in
X^n\}.\]
Note that, since $\{\delta_v\;:\;v\in X^n\}$ is a basis of
$\Bbbk^{X^n}$, the set $\{\delta_v\otimes y_0\;:\;v\in
X^n\}\cup\{\delta_v\otimes y_1\;:\;v\in X^n\}$ is a basis of
$\Bbbk^{X^{n+1}}$. But $\{\delta_v\otimes y_0\;:\;v\in
X^n\}=\{\delta_v\;:\;v\in X^n\}$ is a basis of $\Bbbk^{X^n}$, since
$y_0=\unit$. Consequently, $\mathsf{Y}_{n+1}$ is a basis of $\Bbbk^{X^{n+1}}$.
(Here everywhere $\Bbbk^{X^n}$ are identified with the
corresponding subspaces of $C(X^\omega, \Bbbk)$.)

In the case $\Bbbk=\mathbb{C}$, and identification of $C(X^\omega,
\mathbb{C})$ with a linear subspace of $L^2(X^\omega, \mu)$, where
$\mu$ is the uniform Bernoulli measure on $X^\omega$, it makes sense
to normalize the elements of $\mathsf{Y}_n$ in order to make them of
norm one. Since norm of $\delta_v$ is equal to $2^{-n/2}$, the
recurrent definition of the basis in this case is
\[\mathsf{Y}_{n+1}=\mathsf{Y}_n\cup\{2^{-n/2}\otimes y_1\;:\;v\in
X^n\}.\]
It is easy to see that the union $\mathsf{Y}_\infty$ of the bases
$\mathsf{Y}_n$ is an orthonormal basis of $L^2(X^\omega, \mu)$. It is
called the \emph{Haar basis}. See its use in the context of groups
acting on rooted trees in~\cite{bgr:spec}.
\end{examp}

\section{Automata}

\subsection{Mealy and Moore automata}
\label{ss:MMautomata}
\begin{defi}
\label{def:Mealy}
A \emph{Mealy automaton} (or Mealy machine)
is a tuple \[\mathfrak{A}=(Q, X, Y, \pi, \tau, q_0),\] where
\begin{itemize}
\item $Q$ is the set of \emph{states} of the automaton;
\item $X$ and $Y$ are the \emph{input} and \emph{output alphabets} of
 the automaton;
\item
$\pi:Q\times X\arr Q$ is the \emph{transition map};
\item $\tau:Q\times X\arr X$ is the \emph{output map};
\item $q_0\in Q$ is the \emph{initial state}.
\end{itemize} We always assume that $X$ and $Y$ are finite and have
more than one element each.
\end{defi}

We frequently assume that $X=Y$, and say that the automaton is
\emph{defined over the alphabet $X$}. The automaton is \emph{finite} if
the set $Q$ is finite. In some cases, we do not assume that an initial
state is chosen.

Let $\mathfrak{A}=(Q, X, Y, \pi, \tau, q_0)$ be a Mealy automaton. Let us
extend the definition of the maps $\pi$ and $\tau$ to maps
$\pi:Q\times X^*\arr Q$ and $\tau:Q\times X^*\arr X$ by
the inductive rules
\[\pi(q, xv)=\pi(\pi(q, x), v),\qquad\tau(q, xv)=\tau(\pi(q, x),
v)).\]
We interpret the automaton $\mathfrak{A}$ as a machine, which being in
a state $q\in Q$ and reading a letter $x\in X$, goes to the state
$\pi(q, x)$, and gives the letter $\tau(q, x)$ on the output. If the
machine starts at the state $q\in Q$, and reads a word $v$, then its
final state will be $\pi(q, v)$ is, it the final letter on output will
be $\tau(q, v)$.

\begin{defi}
\label{def:Mealyaction}
The transformation $\mathfrak{A}_{q_0}:X^*\arr X^*$ or
$\mathfrak{A}_{q_0}:X^\omega\arr X^\omega$ defined by a Mealy automaton
$\mathfrak{A}=(Q, X, Y, \pi, \tau, q_0)$ is the map
\begin{equation}\label{eq:Mealyaction}
\mathfrak{A}_{q_0}(x_0x_1x_2\ldots)=\tau(q_0, x_0)\tau(q_1, x_1)\tau(q_2, x_2)\ldots,
\end{equation}
where $q_{i+1}=\pi(q_i, x_i)$.
\end{defi}

In other words, $\mathfrak{A}_{q_0}(v)$ is the word that the machine
gives on output, when it reads the word $v$ on input, if $q_0$ is its
initial state.

\begin{examp}
Let $G\le\autxs$ be a self-similar group. Consider the corresponding
\emph{full automaton} with the set of states $Q=G$, and output and
transition functions defined by the rules:
\[\pi(g, x)=g|_x,\qquad\tau(g, x)=g(x).\]
It follows from~\eqref{eq:sections} that if we choose $g\in G$ as the initial
state, then the transformations of $X^*$ and $X^\omega$ defined by
this automaton coincides with the original transformations defined by
$g\in\autxs$.

This automaton is infinite, but if $G\le\fautxs$, then for every $g\in
G$, the set $\{g|_v\;:\;v\in X^*\}$ is a finite set, and we can take
it as a set of states of a finite automaton defining the
transformation $g$.
\end{examp}

A special type of Mealy automata are the \emph{Moore automata}. The
definition of a Moore automaton is the same as
Definition~\ref{def:Mealy}, except that the output function is a map
$\tau:Q\arr X$, i.e., the output depends only on the state, and does
not depend on the input letter.

Moore automata also act on words, essentially in the same way as Mealy
automata. We can extend the definition of the transition function
$\pi$ to $Q\times X^*$ by the same formula as for the Mealy
automata. Then the action of a Moore automaton with initial state
$q_0$ on words is given by the rule
\begin{equation}
\label{eq:Mooreaction}
\mathfrak{A}_{q_0}(x_1x_2\ldots)=\tau(q_1)\tau(q_2)\ldots,
\end{equation}
where $q_{i+1}=\pi(q_i, x_{i+1})$.

Even though the definition of a Moore automaton seems to be more
restrictive than the definition of a Mealy automaton, the two notions
are basically equivalent, as any Mealy automaton can be modeled by a
Moore automaton. Hence, the set of
maps defined by finite Mealy automata coincides with the set of maps
defined by finite Moore automata.

Let $\mathfrak{A}=(Q, X, Y, \pi, \tau, q_0)$ be a Mealy
automaton. Consider the Moore automaton $\mathfrak{A}'$ over the input
and output alphabets $X$ and $Y$, respectively, with the set of states $Q\times
X\cup\{p_0\}$, where $p_0$ is an element not belonging to $Q\times X$, and
with the transition and output maps $\pi'$ and $\tau'$
given by the rules
\[\pi'(p_0, x)=(\pi(q_0, x), x),\quad\pi'((q, x_1), x_2))=(\pi(q,
x_2), x_2)),\]
and
\[\tau'(q, x_1)=\tau(q, x_1),\]
where $x, x_1, x_2\in X$. (We define $\tau'(p_0)$ to be any
letter, since it will never appear in the output.) It is easy to check
that the new Moore automaton with the initial state $p_0$ defines the
same maps on $X^*$ and $X^\omega$ as the original Mealy automaton
$\mathfrak{A}$.

Therefore, we will not use Moore automata to define transformations of
the sets of words. They will be used to define automatic sequences and
matrices in Section~\ref{s:automatic}. Traditionally, Mealy automata
are used in theory of groups generated by automata
(see~\cite{grineksu_en}), while Moore automata are used for generation
of sequences (even though the term ``Moore automata'' is not used
in~\cite{allouchshallit:sequences}).

\subsection{Diagrams of automata}
The automata are usually represented as finite labeled graphs (called
\emph{Moore diagrams}). The set
of vertices coincides with the set of states $Q$. For every $q\in Q$
and $x\in X$ there is an arrow from $q$ to $\pi(q, x)$ labeled by
$(x, \tau(q, x))$ in the case of Mealy automata, and just by $x$ in the
case of Moore automata. The initial state is marked, and the states
are marked by the values of $\tau(q)$, if it is a Moore automaton.

Sometimes the arrows of diagrams of Mealy automata are just labeled
by the input letters $x$, and the vertices are labeled by the
corresponding transformation $x\mapsto\tau(q, x)$.

Consider a directed graph with one marked (initial) vertex, in
which the edges are labeled by pairs
$(x, y)\in X^2$. The necessary and sufficient condition for such a
graph to represent a Mealy automaton is that for every vertex $q$ and
every letter $x\in X$ there exists a unique arrow starting at $q$ and
labeled by $(x, y)$ for some $y\in X$. Then an image of a word
$x_1x_2\ldots$ under the action of the automaton is calculated by
finding the unique direct path $e_1, e_2, \ldots$ of arrows
starting at the initial vertex, whose
arrows are labeled by $(x_1, y_1), (x_2, y_2), \ldots$,
respectively. Then $y_1y_2\ldots$ is the image of $x_1x_2\ldots$.

\begin{figure}
\centering
\includegraphics{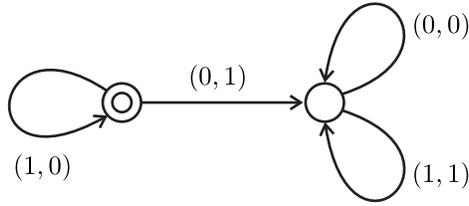}
\caption{The binary adding machine}
\label{fig:admach}
\end{figure}

The diagram of the adding machine transformation (see
Example~\ref{ex:admach}) is shown on Figure~\ref{fig:admach}. We
mark the initial state by a double circle.

\begin{figure}
\centering
\includegraphics{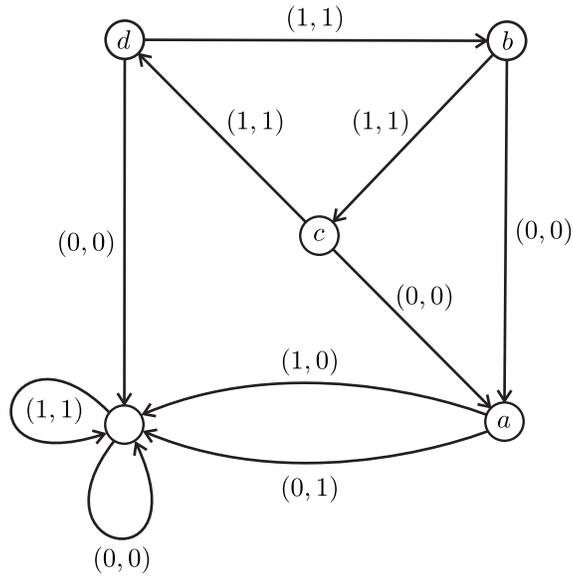}
\caption{The generators of the Grigorchuk group}
\label{fig:grigorch}
\end{figure}

\begin{examp}
The generators $a, b, c, d$ are defined by one automaton, shown on
Figure~\ref{fig:grigorch}, for different choices of the initial
state.
\end{examp}

\subsection{Non-deterministic automata}

Let us generalize the notion of a Mealy automaton by allowing more
general Moore diagrams.

\begin{defi}
A (\emph{non-deterministic}) \emph{synchronous} automaton
$\mathfrak{A}$ over an
alphabet $X$ is an oriented graph whose arrows are labeled by pairs of
letters $(x, y)\in X^2$. Such automaton is called
\emph{$\omega$-deterministic} if for every infinite word
$x_1x_2\ldots\in X^\omega$ and for every vertex (i.e., \emph{state}) $q$
of $\mathfrak{A}$ there exists at most one directed path starting in
$q$ which is labeled by $(x_1, y_1), (x_2, y_2), \ldots$ for some
$y_i\in X$.
\end{defi}

Note that in the above definition for a vertex state $q$ of
$\mathfrak{A}$ and a letter $x\in X$ there maybe several or no edges
starting at $q$ and labeled by $(x, y)$ for $y\in X$. It means that
the automaton $\mathfrak{A}$ may be \emph{non-deterministic} on finite
words and \emph{partial}, i.e., that a state $q$ transforms a finite word $v\in X^*$ into
several different words, and may not accept some of the words on input.

If an automaton $\mathfrak{A}$ is $\omega$-deterministic, then every
its state $q$ defines a map between closed subsets of $X^\omega$,
mapping $x_1x_2\ldots$ to $y_1y_2\ldots$, if there exists a directed
path starting in $q$ and labeled by $(x_1, y_1), (x_2, y_2), \ldots$.

\begin{figure}
\centering
\includegraphics{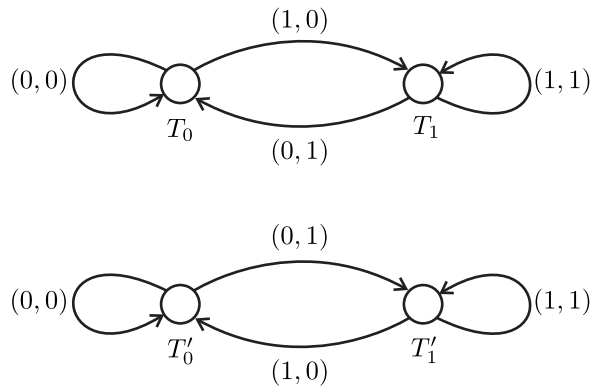}
\caption{Appending and erasing letters}
\label{fig:cuntz}
\end{figure}

\begin{examp}
Let $X=\{0, 1\}$. The states $T_0$ and $T_1$ of the automaton shown on
Figure~\ref{fig:cuntz}
define the transformations $x_1x_2\ldots\mapsto 0x_1x_2\ldots$ and
$x_1x_2\ldots\mapsto 1x_1x_2\ldots$, respectively. The states $T_0'$
and $T_1'$ define the inverse transformations $0x_1x_2\ldots\mapsto
x_1x_2\ldots$ and $1x_1x_2\ldots\mapsto x_1x_2\ldots$.

Note that the first automaton (defining the transformations $T_0$ and
$T_1$) is deterministic. For example, the state $T_0$ acts on the
finite words by transformations $x_1x_2\ldots x_n\mapsto 0x_1x_2\ldots
x_{n-1}$. The second automaton is partial and non-deterministic on
finite words. For example, there are two arrows starting at $T_0'$
labeled by $(0, 1)$ and $(0, 0)$, but no arrows labeled by $(1, y)$.
\end{examp}

An \emph{asynchronous} automaton is defined in the same way, but the
labels are pairs of arbitrary words $(u, v)\in (X^*)^2$.

\begin{defi}
\label{def:synchronautomatichomeo}
A homeomorphism $\phi:X^\omega\arr X^\omega$ is \emph{synchronously}
(resp.\ \emph{asynchronously}) \emph{automatic} if it is defined by a
finite $\omega$-deterministic synchronous (resp.\ asynchronous) automaton.
\end{defi}

A criterion for a homeomorphism to be synchronously automatic is
given in Proposition~\ref{pr:automatichomeo}.

Asynchronously automatic homeomorphisms of $X^\omega$ are studied
in~\cite{grineksu_en,grinek_en}. It is shown there that the set of
synchronously automatic homeomorphisms of $X^\omega$ is a group, and that it
does not depend on $X$ (if $|X|>1$). More precisely, it is proved that
for any two finite alphabets $X, Y$ (such that $|X|, |Y|>1$) there
exists a homeomorphism $X^\omega\arr Y^\omega$ conjugating the
corresponding groups of asynchronously automatic automorphisms. Very
little is known about this group, which is called
in~\cite{grineksu_en} the \emph{group of rational homeomorphisms of
 the Cantor set}.

\section{Automatic matrices}
\label{s:automatic}

\subsection{Automatic sequences}

Here we review the basic definitions and facts about automatic
sequences. More can be found in the monograph~\cite{allouchshallit:sequences,haesler:automatic}.

Let $A$ be a finite alphabet, and let $A^\omega$ be the space of the
right-infinite sequence of elements of $A$ with the direct product
topology.

Fix an integer $d\ge 2$, and consider the transformation
$\stencil_d:A^\omega\arr(A^\omega)^d$, which we called the
\emph{stencil map}, defined by the rule
\[\stencil_d(a_0a_1a_2\ldots)=(a_0a_da_{2d}\ldots,
a_1a_{d+1}a_{2d+1}\ldots, \ldots, a_{d-1}a_{2d-1}a_{3d-1}\ldots).\]
It is easy to see that $\stencil_d$ is a homeomorphism.
We denote the coordinates of $\stencil_d(w)$ by $\stencil_d(w)_i$,
so that
\[\stencil_d(w)=(\stencil_d(w)_0, \stencil_d(w)_1, \ldots,
\stencil_d(w)_{d-1}),\]
and call them \emph{$d$-decimations} of the sequence
$w$. \emph{Repeated $d$-decimations} of $w$ are all sequences that can
be obtained from $w$ by iterative application of the decimation
procedure, i.e., all sequences of the form
\[\stencil_d(\stencil_d(\ldots\stencil_d(w)_{i_n}\ldots)_{i_2})_{i_1}.\]

\begin{defi}
\label{def:automaticsequence}
A sequence $w\in A^\omega$ is \emph{$d$-automatic} if the set of all
repeated $d$-decimations of $w$ (called the \emph{kernel} of $w$
in~\cite[Section~6.6]{allouchshallit:sequences}) is finite.
\end{defi}

We say that a subset $Q\subset A^\omega$ is
\emph{$d$-decimation-closed} if for every $w\in Q$ all
$d$-decimations of $w$ belong to $Q$. The following is
obvious.

\begin{lemma}
A sequence is $d$-automatic if it belongs to a finite
$d$-decimation-closed subset of $A^\omega$.
\end{lemma}

Classically, a sequence $w=a_0a_1\ldots$ is called $d$-automatic if there exists a Moore
automaton $\mathfrak{A}$ with input alphabet $\{0, 1, \ldots, d-1\}$
and output alphabet $A$ such that if $n=i_0+i_1d+\cdots+i_md^m$ is a
base $d$ expansion of $n$, then the output of $\mathfrak{A}$ after
reading the word $i_0i_1\ldots i_m$ is $a_n$. An equivalent variant
of the definition requires that $a_n$ is the output of the automaton
after reading $i_mi_{m-1}\ldots i_1i_0$. One also may allow, or not
$i_m$ to be equal to zero, and the numeration of the letters of the
sequence $w$ to start from 1. All these different definitions of
automaticity of sequences are equivalent to each other,
see~\cite[Section~5.2]{allouchshallit:sequences}. They are also
equivalent to Definition~\ref{def:automaticsequence},
see~\cite[Theorem~6.6.2]{allouchshallit:sequences}.

\begin{figure}
\centering
\includegraphics{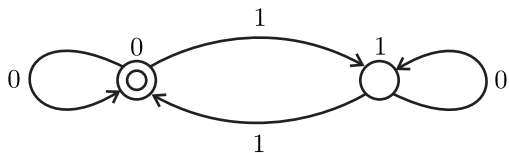}
\caption{Automaton generating the Thue-Morse sequence}
\label{fig:tuemorse}
\end{figure}

\begin{examp}
The \emph{Thue-Morse sequence} is the sequence
$t_0t_1\ldots\in\{0, 1\}^\omega$, where $t_n$ is the sum modulo 2 of
the digits of $n$ in the binary numeration system. The beginning of
length $2^n$ of this sequence can be obtained from $0$ by applying the
substitution \[0\mapsto 01,\qquad 1\mapsto 10\]
$n$ times:
\[0\mapsto 01\mapsto 0110\mapsto 01101001\mapsto
0110100110010110\mapsto\ldots\]
It is easy to see that this sequence is generated by the automaton
shown on Figure~\ref{fig:tuemorse}. Here we label the vertices (the
states) of the automaton by the corresponding values of the output
function. The initial state is marked by a double circle. For more on
properties of the Thue-Morse sequence,
see~\cite[5.1]{allouchshallit:sequences}.
\end{examp}

The last example can be naturally generalized to include all automatic
sequences. Namely, a \emph{$k$-uniform morphism} $\phi:X^*\arr Y^*$
is a morphism of monoids such that $|\phi(x)|=k$ for every $x\in
X$. By a theorem of Combham
(see~\cite[Theorem~6.3.2]{allouchshallit:sequences}) a sequence is
$k$-automatic if and only if it is an image, under a coding (i.e., a
1-uniform morphism), of a
fixed point of a $k$-uniform endomorphism $\phi:X^*\arr X^*$.

\begin{examp}
Consider the alphabet $X=\{a, b, c\}$, and the morphism $\phi:X^*\arr
X^*$ given by
\[\phi(a)=aca,\quad\phi(b)=d,\quad\phi(c)=b,\quad\phi(d)=c.\]
This substitution appears in the presentation~\cite{lysionok:pres} of the Grigorchuk
group.

The fixed point of $\phi$ is obtained as the limit of $\phi^n(a)$, and
starts with $acabacadacabaca\ldots$. The morphism $\phi$ is not
uniform, but it is easy to see that the fixed point belongs to $\{ab,
ac, ad\}^\infty$, and on the words $B=ab, C=ac, D=ad$ it acts on $\{B,
C, D\}^\infty$ as a 2-uniform endomorphism:
\[\phi(B)=acad=CD,\quad\phi(C)=acab=CB,\quad\phi(D)=acac=CC.\] It
follows from Combham's theorem that the fixed point of $\phi$ is
2-automatic.
\end{examp}

Let us show how to construct an automaton producing a sequence
satisfying the conditions of Definition~\ref{def:automaticsequence}.

Suppose that $w_0\in A^\omega$ is automatic, and let $Q$ be
a finite $d$-decimation-closed subset of $A^\omega$ that contains
$w_0$ (for example, we can take $Q$ to be equal to the set
of all repeated $d$-decimations of $w_0$).

Consider a Moore automaton with the set of states $Q$, initial
state $w_0$, input alphabet $\{0, 1, \ldots, d-1\}$, output alphabet
$A$, transition function
\[\pi(w, i)=\stencil_d(w)_i,\]
 output function
\[\tau(x_0x_1\ldots)=x_0.\]
We call the constructed Moore automaton
\[\mathfrak{A}=(Q, \{0, 1, \ldots, d-1\}, A, \pi, \tau, w_0)\]
the \emph{automaton of $w_0$}.

\begin{proposition}
\label{pr:automaticsequences}
Let $w_0=a_0a_1\ldots\in A^\omega$ be an automatic sequence, and let $\mathfrak{A}$
be its automaton. Let $n$ be a non-negative integer, and let $i_0,
i_1, \ldots, i_m$ be a sequence of elements of the set $\{0, 1,
\ldots, d-1\}$ such that $n=i_0+i_1d+i_2d^2+\cdots+i_md^m$. Then
$\tau(w_0, i_0i_1\ldots i_m)=a_n$, i.e., the output of $\mathfrak{A}$
after reading $i_0i_1\ldots i_m$ is $a_n$.
\end{proposition}

\begin{proof}
It follows from the definition of the automaton $\mathfrak{A}$ that
\begin{equation}
\label{eq:piwio1}
\pi(w, i_0\ldots
i_m)=\stencil_d(\ldots\stencil_d(\stencil_d(x_0x_1\ldots)_{i_0})_{i_1}\ldots)_{i_m}
\end{equation}
for all $w=x_0x_1\ldots\in Q$.

It also follows from the definition of the stencil map that the
sequence~\eqref{eq:piwio1} is equal to
\[x_nx_{n+d^{m+1}}x_{n+2d^{m+1 }}x_{n+3d^{m+1}}\ldots,\]
where $n=i_0+i_1d+i_2d^2+\cdots+i_md^m$. It follows that $\tau(w_0,
i_0i_1\ldots i_m)=x_n$.
\end{proof}

\subsection{Automatic infinite matrices}

The notion of automaticity of sequences can be generalized to
matrices in a straightforward way
(see~\cite[Chapter~14]{allouchshallit:sequences},
where they are called
\emph{two-dimensional sequences}).

Let $A$ be a finite alphabet, and let $A^{\omega\times\omega}$ be the
space of all infinite to right and down two-dimensional matrices of elements of
$A$, i.e., arrays of the form
\begin{equation}\label{eq:a}
a=\left(\begin{array}{ccc}a_{11} & a_{12} & \ldots \\ a_{21} & a_{22}
  & \ldots\\ \vdots & \vdots & \ddots\end{array}\right).\end{equation}

\begin{figure}
\centering
\includegraphics[width=3in]{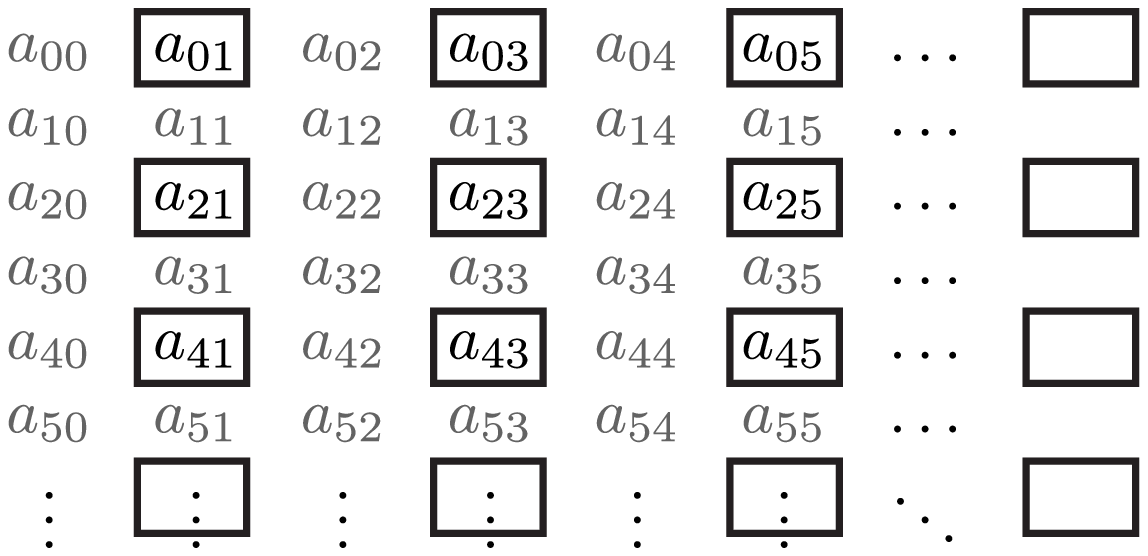}
\caption{The stencil}
\label{fig:stencil}
\end{figure}

Fix an integer $d\ge 2$, and consider the map:
\[\stencil_d:A^{\omega\times\omega}\arr(A^{\omega\times\omega})^{d\times
 d}\]
from the set of infinite matrices to the set of $d\times d$ matrices
whose entries are elements of $A^{\omega\times\omega}$. It maps the
matrix $a$ to the
matrix \[\stencil(a)=\left(\begin{array}{cccc}\stencil_d(a)_{00} &
  \stencil_d(a)_{01} & \ldots & \stencil_d(a)_{0, d-1}\\
  \stencil_d(a)_{10} & \stencil_d(a)_{11} & \ldots &
  \stencil_d(a)_{1, d-1}\\ \vdots & \vdots & \ddots & \vdots \\
  \stencil_d(a)_{d-1, 0} & \stencil_d(a)_{d-1, 1} & \ldots &
  \stencil_d(a)_{d-1, d-1}\end{array}\right),\]
where
\[\stencil_d(a)_{ij}=\left(\begin{array}{cccc}a_{i, j} & a_{i, j+d} &
  a_{i, j+2d} & \ldots\\
a_{i+d, j} & a_{i+d, j+d} & a_{i+d, j+2d} & \ldots\\
a_{i+2d, j} & a_{i+2d, j+d} & a_{i+2d, j+2d} & \ldots\\
\vdots & \vdots & \vdots & \ddots\end{array}\right).\]
The entries of $\stencil_d(a)$ are called $d$-decimations of $a$. We
call $\stencil_d$ the \emph{stencil map}, since entries of the matrix
$\stencil_d(A)$ are obtained from the matrix $A$ by selecting entries
using a ``stencil'' consisting of a square grid of holes, see
Figure~\ref{fig:stencil}  .

The definition of automaticity for matrices is then the same as for
sequences.

\begin{defi}
A matrix $a\in A^{\omega\times\omega}$ is \emph{$d$-automatic} ($[d,
d]$-automatic in terminology of~\cite{allouchshallit:sequences}) if the
set of matrices that can be obtained from $a$ by repeated
$d$-decimations is finite.
\end{defi}

One can also use stencils with a rectangular grid of holes, i.e.,
selecting the entries of a decimation with one step horizontally, and
with a different step vertically. This will lead us to the notion of a
$[d_1, d_2]$-automatic matrix, as in~\cite{allouchshallit:sequences},
but we do not use this notion in our paper.

An interpretation of automaticity of matrices via automata theory is
also very similar to the interpretation for sequences. The only
difference is that the input alphabet of the automaton is the direct
product $\{0, 1, \ldots, d-1\}\times\{0, 1, \ldots, d-1\}$. If we want
to find an entry $a_{n_1, n_2}$ of an automatic matrix defined by a
Moore automaton $\mathfrak{A}$, then we represent the indices $n_1$
and $n_2$ in base $d$:
\[n_1=i_0+i_1d+\cdots+i_md^m,\qquad n_2=j_0+j_1d+\cdots+j_md^m\]
for $0\le i_s, j_t\le d-1$, and then feed the sequence $(i_0,
j_0)(i_1, j_1)\ldots (i_m, j_m)$ to the automaton $\mathfrak{A}$. Its
final output will be $a_{n_1, n_2}$.

We say that a matrix $a=(a_{ij})_{i\ge
 0, j\ge 0}$ over a field $\Bbbk$ is \emph{column-finite} if the number of
non-zero entries in each column of $a$ is finite. The set
$\mat_{\infty}(\Bbbk)$ of all
column-finite matrices is an algebra isomorphic to the algebra of
endomorphisms of the infinite-dimensional vector space
$\Bbbk^\infty=\bigoplus_{\mathbb{N}}\Bbbk$. We denote by $\mat_k(\Bbbk)$ the algebra of $k\times
k$-matrices over $\Bbbk$.

\begin{lemma}
The stencil map
$\stencil_d:\mat_{\infty}(\Bbbk)\arr\mat_d(\mat_{\infty}(\Bbbk))$ is
an isomorphism of $\Bbbk$-algebras.
\end{lemma}

\begin{proof}
A direct corollary of the multiplication rule for matrices.
\end{proof}

Denote by $\mathcal{A}_d(\Bbbk)\subset\mat_\infty(\Bbbk)$ the set of all column-finite
$d$-automatic matrices over $\Bbbk$.

\begin{proposition} Let $\Bbbk$ be a finite field.
Then $\mathcal{A}_d(\Bbbk)$ is an algebra. The stencil
map \[\stencil_d:\mathcal{A}_d(\Bbbk)\arr\mat_d(\mathcal{A}_d(\Bbbk))\]
is an isomorphism of $\Bbbk$-algebras.
\end{proposition}

\begin{proof}
Let $A$ and $B$ be $d$-automatic column-finite matrices. Let
$\mathfrak{A}$ and $\mathfrak{B}$ be finite decimation-closed sets
containing $A$ and $B$, respectively. Then the set
$\{a_1A_1+b_1B_1\;:\;A_1\in\mathfrak{A}, B_1\in\mathfrak{B}, a_1,
b_1\in\Bbbk\}$ is decimation-closed, and it contains any linear
combination of $A$ and $B$. The set is finite, which shows that any
linear combination of $A$ and $B$ is automatic.

Let $\mathfrak{C}$ be the linear span of all products $A_1B_1$ for
$A_1\in\mathfrak{A}$, $B_1\in\mathfrak{B}$. Since the stencil map
$\stencil_d:\mat_\infty(\Bbbk)\arr\mat_d\left(\mat_\infty(\Bbbk)\right)$
is an isomorphism of $\Bbbk$-algebras, the set $\mathfrak{C}$ is
decimation-closed. It contains $AB$, hence $AB$ is automatic.

The map
$\stencil_d:\mathcal{A}_d(\Bbbk)\arr\mat_d\left(\mathcal{A}_d(\Bbbk)\right)$
is obviously a bijection. It is a homomorphism of $\Bbbk$-algebras,
because it is a homomorphism on $\mat_\infty(\Bbbk)$.
\end{proof}

If we have a finite $d$-decimation-closed set of matrices
$\mathfrak{A}$, then its elements are uniquely determined by the
corresponding matrix recursion
$\stencil_d:\mathfrak{A}\arr\mat_d\left(\mathfrak{A}\right)$ and by
the top left entries $a_{00}$ of each matrix
$A\in\mathfrak{A}$. Namely, suppose that we want to find an entry
$a_{m, n}$ of a matrix $A\in\mathfrak{A}$. Let $0\le i<d$ and $0\le
j<d$ be the remainders of division of $n$ and $m$ by $d$. Then $a_{m,
 n}$ is equal to the entry $b_{\frac{m-i}d, \frac{n-j}d}$ of the
matrix $B=\stencil_d(A)_{i, j}$. Repeating this procedure several
times, we eventually will find a matrix $C=(c_{i, j})_{i, j=0}^\infty\in\mathfrak{A}$ such that
$a_{m, n}=c_{0, 0}$.

\begin{examp}
A particular case of automatic matrices are triangular matrices
of the form $$\left(\begin{array}{cccc}a_{00} & a_{01} & a_{02} &
  \ldots\\
0 & a_{11} & a_{12} & \ldots\\
0 & 0 & a_{22} & \ldots\\
\vdots & \vdots & \vdots & \ddots\end{array}\right),$$
where the diagonals $a_{0, i}, a_{1, i+1}, a_{2, i+2}, \ldots$ are
eventually periodic, and only a finite number of them are
non-zero. Note that the set of such uni-triangular matrices is a
group.

The following subgroup of this group (of matrices over the field
$\mathbb{F}_2$) was considered in~\cite{olisush_en}. Let
$B_1=\left(\begin{array}{ccc} 1 & 1 & 1\\ 1 & 0 & 0\\ 1 & 1 &
  1\end{array}\right)$, $C_1=\left(\begin{array}{ccc} 1 & 0 & 0\\ 0
  & 0 & 0\\ 1 & 0 & 0\end{array}\right)$,
$B_2=\left(\begin{array}{ccc} 1 & 1 & 0\\ 0 & 1 & 0\\ 1 & 1 &
  0\end{array}\right)$, $C_2=\left(\begin{array}{ccc} 1 & 0 & 0\\ 1
  & 0 & 0\\ 1 & 0 & 0\end{array}\right)$. It is shown
in~\cite{olisush_en} that the infinite matrices
\[F_1=\left(\begin{array}{cccccc} I & B_1 & C_1 & O & O & \ldots\\
O & I & B_1 & C_1 & O & \ldots\\ O & O & I & B_1 & C_1 & \ldots\\
\vdots & \vdots & \vdots & \vdots & \vdots &
\ddots\end{array}\right),\quad
F_2=\left(\begin{array}{cccccc} I & B_2 & C_2 & O & O & \ldots\\
O & I & B_2 & C_2 & O & \ldots\\ O & O & I & B_2 & C_2 & \ldots\\
\vdots & \vdots & \vdots & \vdots & \vdots &
\ddots\end{array}\right)\]
generate a free group of rank 2. Here $I$ and $O$ are the identity and
zero matrices of size $3\times 3$, respectively.

In Section~\ref{s:unitriangular} we will show that any residually finite $p$-group
can be represented in triangular form. Groups generated by finite
automata will be represented by automatic uni-triangular matrices. The
next subsection is the first step in this direction.
\end{examp}

\subsection{Representation of automata groups by automatic matrices}

Let $\mathsf{B}$ be a basis of $\Bbbk^X$ such that $\unit\in\mathsf{B}$.
Order the elements of $\mathsf{B}$ into a sequence $y_0<y_1<\ldots<y_{d-1}$,
where $y_0=u$. Recall that the inductive limit $\mathsf{B}_\infty$ of the bases
$\mathsf{B}^{\otimes n}$ of $\Bbbk^{X^n}$ with respect to the embeddings $f\mapsto
f\otimes u$ is a basis of $C(X^\omega, \Bbbk)$, whose
elements are infinite tensor products $y_{i_0}\otimes
y_{i_1}\otimes\cdots$, where all but a finite number of factors
$y_{i_k}$ are equal to $y_0=u$. In other words, $\mathsf{B}_\infty$ consists of
functions of the form
\[f(x_0x_1\ldots)=y_{i_0}(x_0)y_{i_1}(x_1)\cdots,\]
where all but a finite number of factors on the right-hand side are
equal to the constant one function.

We can order such products using
the \emph{inverse} lexicographic order, namely $y_{i_0}\otimes
y_{i_1}\otimes\cdots<y_{j_0}\otimes y_{j_1}\otimes\cdots$ if and only
if $y_{i_k}<y_{j_k}$, where $k$ is the \emph{largest} index such that
$y_{i_k}\ne y_{j_k}$.

It is easy to see that the ordinal type of
$\mathsf{B}_\infty$ is $\omega$. Let $e_0<e_1<e_2<\ldots$ be all elements of
$\mathsf{B}_\infty$ taken in the defined order. It is checked directly that if
$e_n=y_{i_0}\otimes y_{i_1}\otimes\cdots$, then $n=i_0+i_1\cdot
d+i_2\cdot d^2+\cdots$, i.e., $\overline{\ldots i_2i_1i_0}$ is the
base-$d$ expansion of $n$ (only a finite number of coefficients $i_j$
are different from zero).

\begin{defi}
\label{def:ssbasis}
An ordered basis $\mathsf{B}$ of $\Bbbk^X$ is called \emph{marked} if
its minimal element is $\unit$.

The \emph{self-similar basis} $\mathsf{B}_\infty$ of $C(X^\omega,
\Bbbk)$ \emph{associated with $\mathsf{B}$} is the inverse
lexicographically ordered set of
functions of the form $e_1\otimes e_2\otimes\ldots$, where
$e_i\in\mathsf{B}$ and all but a finite set of elements $e_i$ are
equal to $\unit$, as it is described above.
\end{defi}

Let $\mathsf{B}=\{b_0, b_1, \ldots, b_{d-1}\}$ be an arbitrary (non
necessarily marked) basis of the space $\Bbbk^X$. We define the
associated \emph{matrix recursion} for linear operators on $C(X^\omega,
\Bbbk)$ in the usual way: given an operator $a$, define its image
$\stencil_{\mathsf{B}}(a)=(A_{i, j})_{i, j=0}^{d-1}$ in $\mat_d(C(X^\omega, \Bbbk)$ by the
rule
\[a(b_j\otimes f)=\sum_{i=0}^{d-1}b_i\otimes A_{i, j}(f),\qquad f\in
C(X^\omega, \Bbbk), 0\le j\le d-1.\]
If $\mathsf{B}$ is the basis $\{\delta_x\}_{x\in X}$, then the
matrix recursion $\stencil_{\mathsf{B}}$ restricted to a self-similar group
$G\le\autxs$ coincides with the matrix recursion~\eqref{eq:matrecursion} coming
directly from the wreath recursion.

\begin{lemma}
Let $\mathsf{B}$ be a marked basis of $\Bbbk^X$.
Then the matrix recursion $\stencil_{\mathsf{B}}$ coincides with the stencil map
for the matrices of linear operators in the associated basis
$\mathsf{B}_\infty$.
\end{lemma}

\begin{proof}
Let $A_{ij}$, $0\le i, j\le d-1$, be the entries of $\Xi_{\mathsf{B}}(a)$. Let
$n=i_0+i_1\cdot d+i_2\cdot d^2+\cdots$ be a non-negative
integer written in base $d$. Then
\begin{equation}
\label{eq:abA}
a(b_{j+dn})=a(b_j\otimes b_n)=\sum_{i=0}^{d-1}b_i\otimes A_{i,
 j}(b_n).
\end{equation}
Let $a_{m, n}$, $0\le m, n<\infty$, be the entries of the matrix of $a$
in the basis $\mathsf{B}_\infty$. Then
\[a(b_{j+dn})=\sum_{k=0}^\infty a_{k,
 j+dn}b_k=\sum_{i=0}^{d-1}\sum_{r=0}^\infty a_{i+dr, j+dn}b_i\otimes
b_r=
\sum_{i=0}^{d-1}b_i\otimes \left(\sum_{r=0}^\infty a_{i+dr, j+dn}b_r\right),\]
which together with~\eqref{eq:abA} implies that $c_{r, n}=a_{i+dr,
 j+dn}$ are the entries of $A_{i, j}$ in the basis
$\mathsf{B}_\infty$, i.e., that $A_{i, j}=\stencil_d(a)_{i, j}$.
\end{proof}

\begin{defi}
Let $\Bbbk$ be a finite field. We say that a linear operator $a$ on
$C(X^\omega, \Bbbk)$ is \emph{automatic} if there exists a finite set
$\mathfrak{A}$ of operators such that $a\in\mathfrak{A}$, and
for every $a'\in\mathfrak{A}$ all entries of the matrix
$\stencil_{\mathsf{B}}(a')$ belong to $\mathfrak{A}$.
\end{defi}

\begin{proposition}
\label{pr:automaticoperators}
Let $\mathsf{B}_1, \mathsf{B}_2$ be two bases of $\Bbbk^X$. An
linear operator on $C(X^\omega, \Bbbk)$ is
automatic with respect to $\mathsf{B}_1$ if and only if it is
automatic with respect to $\mathsf{B}_2$.
\end{proposition}

\begin{proof}
Let $a_1$ be an operator which is automatic with respect to
$\mathsf{B}_1$. Let $\mathfrak{A}$ be the corresponding finite set of
operators, closed with respect to taking entries of the matrix
recursion. Let $\mathfrak{A}'$ be the set of all linear combinations
of elements of $\mathfrak{A}$, which is finite, since we assume that
$\Bbbk$ is finite. If $T$ is the transition matrix from
$\mathfrak{B}_1$ to $\mathfrak{B}_2$, then
\[\stencil_{\mathfrak{B}_2}(a)=T^{-1}\stencil_{\mathfrak{B}_1}(a)T\]
for every linear operator $a$. It follows that $\mathfrak{A}'$ is closed with
respect to taking entries of $\stencil_{\mathfrak{B}_2}$. The set
$\mathfrak{A}'$ is finite, $a_1\in\mathfrak{A}'$, hence $a_1$ is also
automatic with respect to $\mathsf{B}_2$.
\end{proof}

As a direct corollary of Proposition~\ref{pr:automaticoperators} we
get the following relation between finite-state automorphisms of
the rooted tree $X^*$ and automatic matrices.

\begin{theorem}
\label{th:finitestateautomaticmat}
Suppose that $\Bbbk$ is finite. Let $\mathsf{B}$ be a marked basis
of $\Bbbk^X$. Then the
matrix of $\pi_\infty(g)$ in the associated basis $\mathsf{B}_\infty$, where $g\in\autxs$,
is $d$-automatic if and only if $g$ is finite-state.
\end{theorem}

We get, therefore, a subgroup of the group of units of
$\mathcal{A}_d(\Bbbk)$ isomorphic to the group $\fautxs$ of
finite-state automorphisms of the tree $X^*$.

Matrix recursions (i.e., homomorphisms from an algebra $A$ to the
algebra of matrices over $A$) associated with groups acting on rooted
trees, and in more general cases were studied in many papers, for
instance
in~\cite{sid:ring,bgr:spec,bartholdi:ring,sidki:recurmatrices,bacher:determinants,nek:bim,nek:cpalg}.
Note that the algebra generated by the natural representation on
$C(X^\omega, \Bbbk)$ of a group $G$ acting on the rooted tree $X^*$ is
different from the group ring. This algebra (and its analogs) were
studied in~\cite{sid:ring,bartholdi:ring,nek:bim}.

\subsection{Creation and annihilation operators}
\label{ss:creationannihilation}

For $h\in\Bbbk^X$, denote by
$T_h$ the operator on $C(X^\omega, \Bbbk)$ acting by the rule
\[T_h(f)=h\otimes f.\]
It is easy to see that $T_h$ is linear, and that
we have $T_{a_1h_1+a_2h_2}=a_1T_{h_1}+a_2T_{h_2}$ for all $h_1,
h_2\in\Bbbk^X$ and $a_1, a_2\in\Bbbk$.

Consider the dual vector space $(\Bbbk^X)'$ to the space $\Bbbk^X$ of
functions.
We will denote the value of a functional $v$ on a function
$f\in\Bbbk^X$ by $\langle v| f\rangle$. Then for every $v\in
(\Bbbk^X)'$ we have an operator $T_v$ on $C(X^\omega, \Bbbk)$
defined by
\[T_v(f)(x_1, x_2, \ldots)=\langle v| f(x, x_1, x_2, \ldots)\rangle,\]
where $f(x, x_1, x_2, \ldots)$ on the right-hand side of the equation
is seen for every choice of the variables $x_1, x_2, \ldots$ as a
function of one variable $x$, i.e., an element of $\Bbbk^X$.

Let $\mathsf{B}=\{e_i\}_{i=0}^{d-1}$ be a basis of $\Bbbk^X$. Let
$\mathsf{B}'=\{e_i'\}_{i=0}^{d-1}$ be the basis of the dual
space defined by $\langle e_i'|e_j\rangle=\delta_{i, j}$. (Here and
in the sequel, $\delta_{i, j}$ is the Kronecker's symbol equal to 1 when $i=j$, and to 0
otherwise.) We will denote $T_{e_i'}=T_{e_i}'$.

Then $T_{e_i}$ is an isomorphism between the space $C(X^\omega,
\Bbbk)$ and its subspace $e_i\otimes C(X^\omega, \Bbbk)$. It is easy
to see that $T_{e_i}'$ restricted onto $e_i\otimes C(X^\omega, \Bbbk)$
is the inverse of this isomorphism, and that $T_{e_i}'$ restricted onto $e_j\otimes
C(X^\omega, \Bbbk)$ is equal to zero for $j\ne i$.

The operators $T_{e_i}$ and $T_{e_j}'$ satisfy the relations:
\begin{equation}
\label{eq:cuntzrel}
T_{e_i}'T_{e_j}=\delta_{e_i,
 e_j},\qquad\sum_{e_i\in\mathsf{B}}T_{e_i}T_{e_i}'=1.
\end{equation}
The products $T_{e_i}T_{e_i}'$ are projections onto the summands
$e_i\otimes C(X^\omega, \Bbbk)$ of the direct sum decomposition
\[C(X^\omega, \Bbbk)=\bigoplus_{e_i\in\mathsf{B}}e_i\otimes C(X^\omega, \Bbbk).\]

Let $\mathsf{B}$ be a marked basis of $\Bbbk^X$, and let
$\unit=e_0<e_1<\ldots<e_{d-1}$ be its elements. Let
$\mathsf{B}_\infty$ be the associated ordered basis of $C(X^\omega,
\Bbbk)$.

If $A$ is (not necessarily square) finite matrix, then we denote by
$A^{\oplus \infty}$ the infinite matrix
\[\left(\begin{array}{cccc}A & O & O & \ldots\\ O & A & O & \ldots \\
  O & O & A & \ldots\\ \vdots & \vdots & \vdots &
  \ddots \end{array}\right),\]
where $O$ is the zero matrix of the same size as $A$.

The following is a direct corollary of the definitions.

\begin{proposition}
\label{pr:Ei}
The matrices of $T_{e_0}$, $T_{e_1}$, \ldots, $T_{e_{d-1}}$
 in the basis $\mathsf{B}$ are equal to
\[E_0=\left(\begin{array}{c}
1\\
0\\
\vdots\\
0\end{array}\right)^{\oplus\infty},
E_1=\left(\begin{array}{c}
0\\
1\\
\vdots\\
0\end{array}\right)^{\oplus\infty}, \ldots,
E_{d-1}=\left(\begin{array}{cccc}
0\\
0\\
\vdots\\
1\end{array}\right)^{\oplus\infty},
\]
respectively. The matrix of $T_{e_i}'$ is the transpose $E_i^\top$ of the matrix $E_i$.
\end{proposition}

The matrices $E_i$ and $E_i'$ have a natural relation to decimation of
matrices. The proof of the next proposition is a straightforward
computation of matrix products.

\begin{proposition}
\label{pr:recurcuntz}
Let $A=(a_{i, j})_{i, j=0}^\infty$ be an infinite matrix, and let
\[\stencil_d(A)=\left(\begin{array}{cccc}A_{0, 0} & A_{0, 1} & \ldots &
  A_{0, d-1}\\ A_{1, 0} & A_{1, 1} & \ldots & A_{1, d-1}\\
\vdots & \vdots & \ddots & \vdots \\
A_{d-1, 0} & A_{d-1, 1} & \ldots & A_{d-1, d-1}\end{array}\right)\]
be the matrix of its $d$-decimations.
Then
\[A_{i, j}=E_i'AE_j\]
and
\[A=\sum_{i, j=0}^{d-1} E_iA_{i, j}E_j'.\]
\end{proposition}

\begin{corollary}
Let $A$ be an operator on $C(X^\omega, \Bbbk)$, and let
$\mathsf{B}=\{e_i\}$ be a basis of $\Bbbk^X$. Then the entries of the associated matrix
recursion for $A$ are equal to $T_{e_i}'AT_{e_j}$.
\end{corollary}

The next proposition is a direct corollary of Proposition~\ref{pr:Ei}.

\begin{proposition}
\label{pr:Th}
If $h=a_0e_0+a_1e_1+\cdots+a_{d-1}e_{d-1}$, then the matrix of $T_h$
is equal to
$\left(\begin{array}{c}a_0\\ a_1\\ \vdots\\ a_{d-1}\end{array}\right)^{\oplus\infty}.$
\end{proposition}

\begin{corollary}
\label{cor:basisTx}
Let $\mathsf{B}$ be a marked basis of $\Bbbk^X$. Order
the letters of the alphabet $X$ in a sequence $x_0, x_1, \ldots,
x_{d-1}$. Let $T_i=T_{\delta_{x_i}}$ and $T_i'=T_{\delta_{x_i}}'$ be
the corresponding operators, defined using the basis
$\{\delta_{x_i}\}$.

Let $S=(a_{ij})_{i, j=0}^{d-1}$ be the
transition matrix from the basis
$\delta_{x_0}<\delta_{x_1}<\cdots<\delta_{x_{d-1}}$ to $\mathsf{B}$. Let
$S^{-1}=(b_{ij})_{i, j=0}^{d-1}$ be the inverse matrix.

Then the matrix of $T_i$ is \[\left(\begin{array}{c}b_{0, i}\\ b_{1, i}\\
  \ldots \\ b_{d-1, i}\end{array}\right)^{\oplus \infty},\] and the
matrix of $T_i'$ is \[\left(\begin{array}{cccc}a_{i, 0} & a_{i, 1} &
  \ldots & a_{i, d-1}\end{array}\right)^{\oplus \infty}.\]
\end{corollary}

Let us consider now the case of the basis
$\mathsf{D}=\{\delta_x\}_{x\in X}$. For simplicity, let us
denote $T_x=T_{\delta_x}$ and $T_x'=T_{\delta_x}'$. Then the operators
$T_x$ and $T_x'$ act on $C(X^\omega, \Bbbk)$ by the rule
\begin{equation}
\label{eq:cuntzgenerators}
T_x(f)(x_1, x_2, \ldots)=\left\{\begin{array}{ll} f(x_2, x_3, \ldots)
   & \text{if $x=x_1$,}\\ 0 & \text{otherwise.}
\end{array}\right.
\end{equation}
and
\begin{equation}
\label{eq:cuntzgeneratorsdual}
T_x'(f)(x_1, x_2, \ldots)=f(x, x_1, x_2, \ldots).
\end{equation}

In other words, the operator $T_x$ is induced by the natural
homeomorphism $X^\omega\arr xX^\omega:w\mapsto xw$, and $T_x'$ is induced by
its inverse map $xX^\omega\arr X^\omega:xw\mapsto w$.

\begin{proposition}
\label{pr:Txautomatic}
Let $\mathsf{B}_\infty$ be a self-similar basis of $C(X^\omega,
\Bbbk)$ associated with a marked basis $\mathsf{B}$ (see Definition~\ref{def:ssbasis}).
Then the matrices of the operators $T_x$
and $T_x'$, for $x\in X$, in the ordered basis $\mathsf{B}_\infty$ are
$|X|$-automatic.
\end{proposition}

Note that we do not require in this proposition the field $\Bbbk$ to
be finite.

\begin{proof}
Let $\delta_x=\sum_{i=0}^{d-1}\alpha_{x, i}e_i$ for $x\in X$,
$\mathsf{B}=\{e_i\}_{i=0}^{d-1}$ and $\alpha_{x, i}\in\Bbbk$.
It follows from Proposition~\ref{pr:recurcuntz} that the entries of
the matrix
$\stencil_2(T_x)$ with respect to the basis $\mathsf{B}$ are equal to
$\sum_{k=0}^{d-1}\alpha_{x, k}T_{e_i}'T_{e_k}T_{e_j}$. Every product of the form
$T_{e_i}'T_{e_k}T_{e_j}$ is equal, by relations~\eqref{eq:cuntzrel},
either to zero, or to $T_{e_j}$. It follows that decimations of
$T_x$ are either zeros or of the form $\alpha_{x,
  i}T_{e_j}$. It follows that the set of repeated decimations of the
matrix of $T_x$ is contained in $\{\alpha_{x, 0}, \alpha_{x, 1},
\cdots, \alpha_{x, d-1}\}\cdot\{T_{e_0}, T_{e_1}, \ldots, T_{e_{d-1}}\}\cup\{0\}$.
\end{proof}

\begin{examp} Let us consider the case $X=\{0, 1\}$ and
$\mathsf{B}=\{y_0, y_1\}$, where $y_0=\delta_0+\delta_1$ and
$y_1=\delta_1$.  Then the
transition matrix from $\{\delta_0, \delta_1\}$ to $\{y_0, y_1\}$ is
$\left(\begin{array}{cc} 1 & 0\\ 1 & 1\end{array}\right)$, whose
inverse is $\left(\begin{array}{cc} 1 & 0 \\ -1 &
    1\end{array}\right)$.

It follows then from Corollary~\ref{cor:basisTx} that the matrices of
$T_0$, $T_1$, $T_0'$, and $T_1'$ in the basis $\mathsf{B}_\infty$ are
\[T_0=
\left(\begin{array}{r}
   1\\
  -1
\end{array}\right)^{\oplus\infty},\quad
T_1=
\left(\begin{array}{c}
   0 \\
   1 \end{array}\right)^{\oplus\infty},
\]
and
\[
T_0'=\left(\begin{array}{cc}1 &
    0\end{array}\right)^{\oplus\infty},\quad
T_1'=\left(\begin{array}{cc}1 & 1\end{array}\right)^{\oplus\infty}.\]
\end{examp}

\begin{examp}
In the case $\Bbbk=\mathbb{C}$, it is natural to consider the
operators
\[S_x=\frac{1}{\sqrt{|X|}}T_x.\]
Then $S_x$ are isometries of the Hilbert space $L^2(X^\omega)$, and
their conjugates $S_x^*$ are equal to $\sqrt{|X|}T_x'$.

The $C^*$-algebra of operators on $L^2(X^\omega)$ generated by the
operators $S_x$ is called the Cuntz algebra~\cite{cuntz}, and is usually denoted
$\mathcal{O}_{|X|}$. Any isometries
satisfying the relations
\[S_x^*S_x=1,\qquad\sum_{x\in X}S_xS_x^*=1\]
generate a $C^*$-algebra isomorphic to $\mathcal{O}_{|X|}$. In
particular, the $C^*$-algebra generated by the matrices $E_i$ is the
Cuntz algebra. Representation of the Cuntz algebra by matrices $E_i$
is an example of a \emph{permutational representation} of
$\mathcal{O}_d$. More on such and similar representations, see~\cite{cuntz_rep}.

Recall that, for $X=\{0, 1\}$,
the Walsh basis of $L^2(X^\omega)$ is the basis $\mathsf{W}_\infty$
constructed starting from the basis $\mathsf{W}=\{y_0, y_1\}$, where
$y_0=\delta_0+\delta_1$ and $y_1=\delta_0-\delta_1$. Then direct
computation with of the transition matrices show that the matrices of $S_0$ and $S_1$ are
\[S_0=\left(\begin{array}{c}\frac{1}{\sqrt{2}}\\
    \frac{1}{\sqrt{2}}\end{array}\right)^{\oplus \infty},\qquad
S_1=\left(\begin{array}{r}\frac{1}{\sqrt{2}} \\
    -\frac{1}{\sqrt{2}}\end{array}\right)^{\oplus \infty}.\]
\end{examp}

\subsection{Cuntz algebras and Higman-Thompson groups}

If $\psi:X^\omega\arr X^\omega$ is a homeomorphism, then it induces a linear
operator $L_\psi$ on $C(X^\omega, \Bbbk)$ given by
\[L_\psi(f)(w)=f(\psi^{-1}(w))\]
for $f\in C(X^\omega, \Bbbk)$ and $w\in X^\omega$.

Fixing any ordered basis of $C(X^\omega, \Bbbk)$,
we get thus a natural faithful representation of the homeomorphism
group of $X^\omega$ in the group of units of the algebra
$\mat_\infty(\Bbbk)$ of column-finite matrices over $\Bbbk$.

\begin{proposition}
\label{pr:automatichomeo}
Let $\psi$ be a homeomorphism of $X^\omega$. Let $u, v\in X^n$, and
denote by $\psi_{u, v}$ the partially defined map given by the formula
\[\psi_{u, v}(w)=\left\{\begin{array}{ll} w' & \text{if $\psi(vw)=uw'$,}\\
  \text{not defined} & \text{otherwise.}\end{array}\right.\]

The following conditions are equivalent.
\begin{enumerate}
\item The homeomorphism $\psi$ is synchronously automatic.
\item The set of partial maps $\{\psi_{u, v}\;:\;u, v\in X^*,
  |u|=|v|\}$ is finite.
\item For every finite field $\Bbbk$, the operator
  $L_\psi:C(X^\omega, \Bbbk)\arr C(X^\omega, \Bbbk)$ is automatic.
\item For some finite field $\Bbbk$, the operator
 $L_\psi:C(X^\omega, \Bbbk)\arr C(X^\omega, \Bbbk)$ is automatic.
\end{enumerate}
\end{proposition}

Synchronously automatic homeomorphisms are defined in
Definition~\ref{def:synchronautomatichomeo}.

\begin{proof}
Equivalence of conditions (2), (3), and (5) follow directly from
Proposition~\ref{pr:recurcuntz}.

Suppose that $\psi$ is synchronously automatic. Let $\mathfrak{A}$
be an initial automaton defining $\psi$. For every
pair $u=a_1a_2\ldots a_n, v=b_1b_2\ldots b_n\in X^*$ of words of equal
length, let $Q_{u, v}$ be the set
of states $q$ of $\mathfrak{A}$ such that there exists a directed path
starting in the initial state $q_0$ of $\mathfrak{A}$ and labeled by
$(a_1, b_1), (a_2, b_2), \ldots, (a_n, b_n)$. Then the set $Q_{u, v}$
defines the map $\psi_{u, v}$ in the following sense. We have
$\psi_{u, v}(x_1x_2\ldots)=y_1y_2\ldots$ if and only if there exists a
path starting in an element of $Q_{u, v}$ and labeled by $(x_1, y_1),
(x_2, y_2), \ldots$. It follows that the number of possible maps of
the form $\psi_{u, v}$ is not larger than the number of subsets of the
set of states of $\mathfrak{A}$. This shows that every synchronously
automatic homeomorphism satisfies condition (2).

Suppose now that a homeomorphism $\psi$ satisfies condition (2), and
let us show that it is synchronously automatic. Construct an automaton
$\mathfrak{A}$
with the set of states $Q$ equal to the set of non-empty maps of the form
$\psi_{u, v}$. For every $(x, y)\in X^2$ and $\psi_{u, v}\in Q$ we
have an arrow from $\psi_{u, v}$ to $\psi_{ux, vy}$, labeled by
$(x, y)$, provided the map $\psi_{ux, vy}$ is not empty. The initial
state of the automaton is the map $\psi=\psi_{\emptyset,
  \emptyset}$. Let us show that this automaton defines the
homeomorphism $\psi$. It is clear that if
$\psi(x_1x_2\ldots)=y_1y_2\ldots$, then there exists a path starting
at the initial state of $\mathfrak{A}$ and labeled by $(x_1, y_1),
(x_2, y_2), \ldots$. On the other hand, if such a path exists for a
pair of infinite words $x_1x_2\ldots, y_1y_2\ldots$, then the maps
$\psi_{x_1x_2\ldots x_n, y_1y_2\ldots y_n}$ are non-empty for every
$n$. In other words, for every $n$ the set $W_n$ of infinite sequences
$w\in x_1x_2\ldots x_nX^\omega$
such that $\psi(w)\in y_1y_2\ldots y_n X^\omega$ is
non-empty. It is clear that the sets $W_n$ are closed and
$W_{n+1}\subset W_n$ for every $n$. By compactness of $X^\omega$ it
implies that $\bigcap_{n\ge 1}W_n$ is non-empty. It follows that
$\psi(x_1x_2\ldots)=y_1y_2\ldots$.
\end{proof}

The next corollary follows directly from condition (2) of Proposition~\ref{pr:automatichomeo}.

\begin{corollary}
The set of all automatic homeomorphisms of $X^\omega$ is a group.
\end{corollary}

We have already seen in Theorem~\ref{th:finitestateautomaticmat}
that a homeomorphisms $g$ of $X^\omega$ defined by an
automorphisms of $X^*$ is automatic if and only if it is
finite state. Note that in this case
$g_{u, v}$ is either empty (if $g(v)\ne u$) or is equal to $g|_v$.

Another example of a group of finitely automatic homeomorphisms of
$X^\omega$ is the Higman-Thompson group $\mathcal{V}_{|X|}$. It is the set of all
homeomorphisms that can be defined in the following way. We say that a
subset $A\subset X^*$ is a \emph{cross-section} if the sets $uX^\omega$
for $u\in A$ are disjoint and their union is $X^\omega$. Let $A=\{v_1,
v_2, \ldots, v_n\}$ and $B=\{u_1, u_2, \ldots, u_n\}$ be cross-sections
of equal cardinality together with a bijection $v_i\mapsto u_i$.
Define a homeomorphism $\psi:X^\omega\arr
X^\omega$ by the rule
\begin{equation}\label{eq:psi}
\psi(v_iw)=u_iw.
\end{equation}
The set of all homeomorphisms that can be defined in this way is the
Higman-Thompson group $\mathcal{V}_{|X|}$, see~\cite{thompson,intro_tomp}.

Let $\psi$ be the homeomorphism defined by~\eqref{eq:psi}. It follows
directly from~\eqref{eq:cuntzgenerators}
and~\eqref{eq:cuntzgeneratorsdual} that the operator $L_\psi$
 induced by $\psi$ is equal to
\[L_\psi=\sum_{i=1}^n T_{u_i}T_{v_i}',\]
where we use notation
\[T_{x_1x_2\ldots x_m}=T_{x_1}T_{x_2}\cdots T_{x_n},\qquad
T_{x_1x_2\ldots x_m}'=T_{x_m}'T_{x_{m-1}}'\cdots T_{x_1}'.\]

The next proposition follows then from Proposition~\ref{pr:Txautomatic}.

\begin{proposition}
The Higman-Thompson group $\mathcal{V}_{|X|}$ is a subgroup of the group of
synchronously automatic homeomorphisms of $X^\omega$.
\end{proposition}

The group generated by $\mathcal{V}_2$ and the Grigorchuk group was
studied by K.~Roever in~\cite{roever}. He proved that it is a finitely
presented simple group isomorphic to the
abstract commensurizer of the Grigorchuk group. Generalizations of
this group (for arbitrary self-similar group) was studied
in~\cite{nek:bim}.

\section{Representations by uni-triangular matrices}
\label{s:unitriangular}
\subsection{Sylow $p$-subgroup of $\autxs$}
\label{ss:sylow}

Let $|X|=p$ be prime. We assume that $X=\{0, 1, \ldots, p-1\}$ is
equal to the field $\fieldp$ of $p$ elements. From now on, we will write vertices
of the tree $X^*$ as tuples $(x_1, x_2, \ldots, x_n)$ in order not to
confuse them with products of elements of $\fieldp$.

Denote by $\sylow$ the subgroup of $\autxs$ consisting of
automorphisms $g$ whose labels $\alpha_{g, v}$ of the vertices of
the portrait consist
only of powers of the cyclic permutation $\sigma=(0, 1, \ldots, p-1)$.
It follows from~\eqref{eq:tablproduct} and~\eqref{eq:tablinverse} that
$\sylow$ is a group. The study of the group $\sylow$ (and its finite
analogs were initiated by
L.~Kaloujnine~\cite{kalouj:sylowfin,kalouj:pinfty-cras,kaloujnine:pinfty}).

Suppose that an element $g\in\sylow$ is represented by a tableau
\[[a_0, a_1(x_1), a_2(x_1, x_2), \ldots],\] as in
Subsection~\ref{ss:rootedtrees}.
Then $a_n(x_1,
x_2, \ldots, x_n)$ are maps from $X^n$ to the group generated by
the cyclic permutation $\sigma$. The elements of this
group act on $X=\fieldp$ by maps $\sigma^a:x\mapsto x+a$. It follows that we
can identify functions $a_n$ with maps $X^n\arr\fieldp$, so that an
element $g\in\sylow$ represented by a
tableau \[[a_0, a_1(x_1), a_2(x_1, x_2), \ldots]\] acts on sequences
$v=(x_0, x_1, \ldots)\in X^\omega$  by the rule
\[g(v)=(x_0+a_0, x_1+a_1(x_0), x_2+a_2(x_1, x_2), x_3+a_3(x_1, x_2,
x_3), \ldots).\]

It follows that if $g_1, g_2\in\sylow$ are represented by the tableaux
$[a_n]_{n=0}^\infty$ and $[b_n]_{n=0}^\infty$, then their product
$g_1g_2$ is represented by the tableau
\begin{multline}
\label{eq:tablproductsylow}
[b_0+a_0, \quad
b_1(x_0)+a_1(x_0+b_0), \\
b_2(x_0, x_1)+a_2(x_0+b_0, x_1+b_1(x_0)),
\\ b_3(x_0, x_1, x_2)+a_3(x_0+b_0, x_1+b_1(x_0), x_2+b_2(x_0, x_1)), \quad\ldots].
\end{multline}

Denote by $\sylown$ the quotient of $\sylow$ by the pointwise
stabilizer of the $n$th level of the tree $X^*$. We can consider
$\sylown$ as a subgroup of the automorphism group of the finite
subtree $X^{[n]}=\bigcup_{k=0}^nX^k\subset X^*$.

\begin{proposition}
The group $\sylown$ is a Sylow subgroup of  the symmetric group
$\symm[X^n]$ and of the automorphism group of the tree $X^{[n]}$.
\end{proposition}

\begin{proof}
The order of $\symm[X^n]$ is $p^n!$, and the maximal power of $p$
dividing it is
\[\frac{p^n}{p}+\frac{p^n}{p^2}+\cdots+\frac{p^n}{p^n}=\frac{p^n-1}{p-1}.\]
It follows that the order of the Sylow $p$-subgroup of $\symm[X^n]$ is $p^{\frac{p^n-1}{p-1}}$.
The order of $\sylown$ is equal to the number of possible tableaux
\[[a_0, a_1(x_1), a_2(x_1, x_2), \ldots a_{n-1}(x_1, \ldots,
x_{n-1})],\] where $a_i$ is an arbitrary map from $X^i$ to the cyclic
group $\langle\sigma\rangle$ of order $p$. The number of possibly maps $a_i$ is hence
$p^{p^i}$. Consequently, the number of possible tableaux is
$p^{1+p+p^2+\cdots+p^{n-1}}=p^{\frac{p^n-1}{p-1}}$.
Since the group of all automorphisms of the tree $X^{[n]}$ is
contained in $\symm[X^n]$ and contains $\sylown$, the subgroup
$\sylown$ is its Sylow $p$-subgroup too.
\end{proof}

\begin{proposition}
Let $g\in\sylow$ be represented by a tableau $[a_0, a_1(x_1), a_2(x_1,
x_2), \ldots]$. Consider the map $\alpha:\sylow\arr\fieldp^\omega$,
where $\fieldp^\omega$ is the infinite Cartesian product of additive
groups of $\fieldp$, given by
\begin{equation}
\label{eq:alphamap}
\alpha(g)=\left(a_0, \sum_{x_1\in\fieldp}a_1(x_1), \sum_{(x_1,
  x_2)\in\fieldp^2}a_2(x_1, x_2), \ldots\right).
\end{equation}
In other words, we just sum up modulo $p$ all the decorations of the portrait of
$g$ on each level. Then $\alpha$ is the abelianization epimorphism
$\sylow\arr\sylow/[\sylow, \sylow]\cong\fieldp^\omega$.
\end{proposition}

\begin{proof}
It is easy to check that $\alpha$ is a homomorphism. It remains to
show that its kernel is the derived subgroup of $\sylow$. This is a
folklore fact, and we show here how it follows from a more general
result of Kaloujnine.

Let $g\in\sylow$ be represented by a tableau
\[[a_0, a_1(x_1), a_2(x_1, x_2), \ldots].\]
Each function $a_n(x_1, x_2, \ldots, x_n)$ can be written as a
polynomial \[\sum_{0\le k_i\le p-1}c_{k_1, k_2, \ldots,
  k_n}x_1^{k_1}x_2^{k_2}\cdots x_n^{k_n}\] for some coefficients $c_{k_1, k_2,
  \ldots, k_n}\in\fieldp$.

It is proved in~\cite[Theorem 6]{kalouj:sylowfin} (see also Equation
(5,4) in~\cite{kaloujnine:pinfty}) that the derived subgroup $[\sylow,
\sylow]$ of $\sylow$ is the set of elements defined by tableaux in which
$a_0$ and the coefficient $c_{p-1, p-1, \ldots, p-1}$ at the
eldest term $x_1^{p-1}x_2^{p-1}\ldots x_n^{p-1}$ are equal to zero for every $n$.

Note that $\sum_{x\in\fieldp}x^k$ is equal to zero
for $k=0, 1, \ldots, p-2$ and is equal to $-1$ for $k=p-1$. Therefore,
\[\sum_{(x_1, x_2, \ldots,
  x_n)\in\fieldp^n}x_1^{k_1}x_2^{k_2}\cdots
x_n^{k_n}=\prod_{i=1}^n\sum_{x\in\fieldp}x^{k_i}\]
is equal to zero for all $n$-tuples $(k_1, k_2, \ldots, k_n)\in\{0, 1,
\ldots, p-1\}^n$ except for $(p-1, p-1, \ldots, p-1)$, when it is
equal to $(-1)^n$. It follows that the coefficient at the eldest term
of $a_n(x_1, x_2, \ldots, x_n)$ is equal to zero if and only if
$\sum_{(x_1, x_2, \ldots, x_n)\in\fieldp^n}a_n(x_1, x_2, \ldots, x_n)=0$.
\end{proof}

\subsection{Polynomial bases of $C(X^\omega, \fieldp)$}
\label{ss:polynombases}

\begin{proposition}
\label{pr:triangular}
Suppose that an ordered basis $\mathsf{B}$ of $\Bbbk^X$ is such that
the matrices of $\pi_1(g)$ for $g\in G\le\autxs$ are all upper
uni-triangular, and the
minimal element of $\mathsf{B}$ is the constant one function. Then the
matrices of $\pi_n(g)$ in the basis $\mathsf{B}^{\otimes n}$ and of
$\pi_\infty(g)$ in the associated basis $\mathsf{B}_\infty$ are
upper uni-triangular for all $g\in G$.
\end{proposition}

See Subsection~\ref{ss:indlimit} for definition of the representations
$\pi_n$. We say that a matrix is \emph{upper uni-triangular} if all
its elements below the main diagonal are equal to zero, and all
elements on the diagonal are equal to one. From now on, unless the contrary is
specifically mentioned, ``uni-triangular'' will mean ``upper uni-triangular''.

\begin{proof}
Let $b_0<b_1<\ldots<b_{d-1}$ be the ordered basis $\mathsf{B}$.
Let $\mathsf{Y}=\{y_0, y_1, \ldots, y_{d-1}\}$ be the
corresponding basis of the right module
$\bim_G$. Namely, we take for every $\sum_{x\in
  X}a_x\delta_x\in\mathsf{B}$ the corresponding element
$\sum_{x\in X}a_xx\in\bim=\Bbbk[X\cdot G]$, see
Subsection~\ref{ss:ssbimodule}. Then
$\mathsf{B}=\mathsf{Y}\otimes\varepsilon$, where $[\varepsilon]$ is
the  left $G$-module of the trivial representation of $G$, see~\ref{ss:indlimit}.

If the matrices of $\pi_1(g)$ are uni-triangular, then
\[\pi_1(g)(b_i)=b_i+a_{i-1, i}b_{i-1}+a_{i-2, i}b_{i-2}+\cdots+a_{0,
  i}b_0\] for some $a_{k, i}\in\Bbbk$ and all $i$. It follows that, in
the bimodule $\bim$, we have relations
\begin{equation}
\label{eq:rectriangle}
g\cdot y_i=y_i\cdot g_{i, i}+a_{i-1, i}y_{i-1}\cdot g_{i-1,
  i}+a_{i-2, i}y_{i-2}\cdot g_{i-2, i}+\cdots+a_{0, i}y_0\cdot g_{0,
  i}
\end{equation}
for some $g_{j, i}\in\Bbbk[G]$ such that $g_{j,
  i}\cdot\epsilon=\epsilon$. (Recall that the last equality just
means that the sum of coefficients of $g_{j, i}\in\Bbbk[G]$, i.e., the
value of the augmentation map, is equal to one.)
Consequently, relation~\eqref{eq:rectriangle} together with the
condition $g_{j, i}\cdot\varepsilon=\varepsilon$ hold for all
$g\in\Bbbk[G]$ such that $g\cdot\varepsilon=\varepsilon$.

It follows that every element $g\cdot
y_{i_1}\otimes y_{i_2}\otimes\cdots\otimes y_{i_n}\in\bim^{\otimes n}$ is equal to
$y_{i_1}\otimes y_{i_2}\otimes\cdots\otimes y_{i_n}\cdot h$  plus a
sum of elements of the form $y_{j_1}\otimes
y_{j_2}\otimes\cdots y_{j_n}\cdot a_{j_1, j_2, \ldots, j_n}$, where $h\in\Bbbk[G]$ is such
that $h\cdot\varepsilon=\varepsilon$, $a_{j_1, j_2, \ldots,
  j_n}\in\Bbbk[G]$, and $j_k\le i_k$ for all $k=1, 2, \ldots,
n$, and $(j_1, j_2, \ldots, j_n)\ne (i_1, i_2, \ldots, i_n)$.
Taking tensor product with $\varepsilon$ and applying
Proposition~\ref{pr:bimeps}, we conclude that for every function
$b_{i_1}\otimes b_{i_2}\otimes \cdots\otimes b_{i_n}\in\Bbbk^{X^n}$
the function $\pi_\infty(g)(b_{i_1}\otimes b_{i_2}\otimes\cdots\otimes
b_{i_n})$ is equal to $b_{i_1}\otimes b_{i_2}\otimes\cdots\otimes
b_{i_n}$ plus a linear combination of functions $b_{j_1}\otimes b_{j_2}\otimes\cdots\otimes
b_{j_n}$ such that $j_k\le i_k$ for all $k=1, 2, \ldots,
n$, and $(j_1, j_2, \ldots, j_n)\ne (i_1, i_2, \ldots, i_n)$. But any
such function $b_{j_1}\otimes b_{j_2}\otimes\cdots\otimes
b_{j_n}$ is an element of $\mathsf{B}_\infty$, which is smaller than
$b_{i_1}\otimes b_{i_2}\otimes\cdots\otimes b_{i_n}$ in the inverse
lexicographic order. This proves that the matrix of $\pi_\infty(g)$ in
the basis $\mathsf{B}_\infty$ is uni-triangular.
\end{proof}

Throughout the rest of our paper we assume that $|X|=p$ is prime,
$\Bbbk$ is the field $\fieldp$ of $p$ elements, $G$ is a subgroup of $\sylow$, and
we identify $X$ with $\fieldp$. We will be able then to use
Proposition~\ref{pr:triangular} to construct bases of $C(X^\omega,
\fieldp)$ in which the representation $\pi_\infty$ of $\sylow$ (and
hence of $G$) are uni-triangular.

Every function $f\in\fieldp^X$ can be represented as a polynomial
$f(x)\in\fieldp[x]$, using the formula
\[\delta_t(x)=\frac{x(x-1)(x-2)\cdots(x-p+1)}{(x-t)},\]
where $(x-t)$ in the numerator and the denominator cancel each
other. (Recall that $(p-1)!=-1\pmod{p}$, by Wilson's theorem.)

Since $x^p=x$ as a function on $\fieldp$ (by Fermat's little theorem),
representations that differ by an element of the ideal generated by
$x^p-x$ represent the same function. Note that the ring
$\fieldp[x]/(x^p-x)$ has cardinality $p^p$, hence we get a natural
bijection between $\fieldp[x]/(x^p-x)$ and $\fieldp^X$, mapping a
polynomial to the function it defines on $\fieldp$. From now on, we
will thus identify the space of functions $\fieldp^X$ with the
$\fieldp$-algebra $\fieldp[x]/(x^p-x)$.

Following Kaloujnine, we will call the
elements of $\fieldp[x]/(x^p-x)$ \emph{reduced polynomials}. We write
them as usual polynomials $a_0+a_1x+\cdots +a_{p-1}x^{p-1}$ (but
keeping in mind reduction, when performing multiplication).

Suppose that $g\in G$ is such
that $g(x)=x+1$ for all $x\in X$. Then $\pi_1(g)$ acts on the
functions $f\in V_1=\fieldp^X$ by the rule
\[\pi_1(g)(f)(x)=f(x-1).\]

In particular, if we represent $f$ as a polynomial, then $\pi_1(g)$
does not change its degree and the coefficient of the leading term.
It follows that the matrix of the operator $\pi_1(g)$ in
the basis $e_0(x)=\unit, e_1(x)=x, e_2(x)=x^2, \ldots, e_{p-1}(x)=x^{p-1}$ is
uni-triangular. Let us denote this marked basis by $\mathsf{E}$.

\begin{defi}
The basis $\mathsf{E}_\infty$ of $C(X^\omega,
\fieldp)$ corresponding to $\mathsf{E}$ and
consisting of all monomial functions $x_1^{k_1}x_2^{k_2}\cdots$
on $X^\omega$ ordered inverse lexicographically, so that
\[e_0=1, e_1=x_1, e_2=x_1^2, \ldots, e_{p-1}=x_1^{p-1}, e_p=x_2,
e_{p+1}=x_1x_2, \ldots\]
is called the \emph{Kaloujnine basis} of monomials.
\end{defi}

It is easy to see that $e_n\in\mathsf{E}_\infty$ is equal to the monomial
function
\begin{equation}\label{eq:en}
e_n(x_1x_2\ldots)=x_1^{k_1}x_2^{k_2}\cdots,\end{equation}
where $k_1, k_2, \ldots$ are the digits of
the base $p$ expansion of $n$, i.e., such that $k_i\in\{0, 1, \ldots,
p-1\}$ and
\[n=k_1+k_2\cdot p+k_3\cdot p^2+\cdots.\]

Coordinates of a function $f\in C(X^\omega, \fieldp)$ in the basis
$\mathsf{E}_\infty$ are the coefficients of the representation of the
function $f$ as a polynomial in the variables $x_1, x_2,
\ldots$. Since we are dealing with functions, we assume that these
polynomials are reduced, i.e., are elements of the ring $\fieldp[x_1,
x_2, \ldots]/(x_1^p-x_1, x_2^p-x_2, \ldots)$.

As an immediate corollary of Proposition~\ref{pr:triangular} we get
the following.

\begin{theorem}
The representation $\pi_\infty$ of $\sylow$ in the Kaloujnine basis
$\mathsf{E}_\infty$ is uni-triangular. In particular, the
representations $\pi_\infty$ of all its subgroups $G\le\sylow$ are
uni-triangular in $\mathsf{E}_\infty$.
\end{theorem}

We can change the ordered basis $\mathsf{E}=\{e_0=\unit, e_1=x, \ldots,
e_{p-1}(x)=x^{p-1}\}$ to any ordered
basis $\mathsf{F}=(f_0, f_1, \ldots, f_{p-1})$ consisting of polynomials of
degrees $0, 1, 2, \ldots, p-1$, respectively, since then the
transition matrix from $\mathsf{E}$ to $\mathsf{F}$ will be triangular, hence the
representation of $G$ in the basis $\mathsf{F}$ will be also uni-triangular.

For example, a natural choice is the basis $\mathsf{B}$ in which the matrix of
the cyclic permutation $x\mapsto x+1:X\arr X$ is the Jordan cell
\[\left(\begin{array}{cccc} 1 & 1 & 0 & \ldots\\ 0 & 1 & 1 & \ldots\\
    0 & 0 & 1 & \ldots\\ \vdots & \vdots & \vdots &
    \ddots\end{array}\right).\]

To get such a basis, define the functions $b_0, b_1, \ldots, b_{p-1}\in V_1$ by the
formula
\[b_k=(\pi_1(g)-1)^{p-1-k}(\delta_0)\]
for $k=0, 1, \ldots, p-2$ and $b_{p-1}=\delta_0$.
Then $(\pi_1(g)-1)(b_k)=(\pi_1(g)-1)^{p-k}(\delta_0)=b_{k-1}$, for all
$k=1, 2, \ldots, p-1$, i.e.,
\[b_{k-1}(x)=b_k(x-1)-b_k(x).\]

Note that
\[b_0=(\pi_1(g)-1)^{p-1}(\delta_0)=(\pi_1(g)^{p-1}+\pi_1(g)^{p-2}+\cdots+1)(\delta_0)=\sum_{k=0}^{p-1}\delta_k=\unit,\]
i.e., the basis $b_0<b_1<\ldots<b_{p-1}$ is marked.

\begin{proposition}
\label{pr:fk}
For every $k\in\{1, \ldots, p-1\}$ and $x\in X=\fieldp$ we have
\[b_k(x)=(-1)^k\binom{x+k}{k}=(-1)^k\frac{(x+1)(x+2)\cdots (x+k)}{k!}.\]
\end{proposition}

Note that $k!\ne 0$ in $\fieldp$ for every $k=1, 2, \ldots,
p-1$.

\begin{proof}
We have $(p-1)!=1$ and $(-1)^{p-1}=1$ in $\fieldp$. We also have
$(x+1)(x+2)\cdots (x+p-1)=0$ for all $x\in\fieldp\setminus\{0\}$. It
follows that $(-1)^{p-1}\binom{x+p-1}{p-1}=\delta_0=f_{p-1}$.

It is enough now to check that the functions $(-1)^k\binom{x+k}{k}$
satisfy the recurrent relation $b_{k-1}(x)=b_k(x-1)-b_k(x)$.
But we have
\begin{multline*}
(-1)^k\binom{x-1+k}{k}-(-1)^k\binom{x+k}{k}=\\
(-1)^{k-1}\left(\binom{x+k}{k}-\binom{x+k-1}{k}\right)=(-1)^{k-1}\binom{x+k-1}{k-1}
\end{multline*}
by the well known identity
\[\binom{a}{b}=\binom{a-1}{b}+\binom{a-1}{b-1}.\]
\end{proof}

\begin{proposition}
\label{pr:transitionmatrixb}
The transition matrix from the basis $(\delta_0, \delta_1, \ldots,
\delta_{p-1})$ to the basis $\mathsf{B}=(b_0, b_1, \ldots, b_{p-1})$ is
\[T=\left(\begin{array}{cccccc}
\binom{p-1}{0} & \binom{p-1}{1} &
    \binom{p-1}{2} & \binom{p-1}{3} & \ldots & \binom{p-1}{p-1} \\
\binom{p-2}{0} & \binom{p-2}{1} &
    \binom{p-2}{2} & \binom{p-2}{3} & \ldots & 0\\
\vdots & \vdots & \vdots & \vdots & \ddots & \vdots\\
1 & 3 & 3 & 1 & \ldots & 0\\
1 & 2 & 1 & 0 & \ldots & 0\\
1 & 1 & 0 & 0 & \ldots & 0\\
1 & 0 & 0 & 0 & \ldots & 0\end{array}\right)\]
Its inverse is obtained by transposing $T$ with respect to the
secondary diagonal:
\[T^{-1}=\left(\begin{array}{cccccc} 0 & 0 & 0 & 0 & \ldots &
    \binom{p-1}{p-1}\\
0 & 0 & 0 & 0 & \ldots & \binom{p-1}{p-2} \\
\vdots & \vdots & \vdots & \vdots & \ddots & \vdots \\
0 & 0 & 0 & 1 & \ldots & \binom{p-1}{3}\\
0 & 0 & 1 & 3 & \ldots & \binom{p-1}{2}\\
0 & 1 & 2 & 3 & \ldots & \binom{p-1}{1}\\
1 & 1 & 1 & 1 & \ldots & \binom{p-1}{0}\end{array}\right).\]
\end{proposition}

\begin{proof}
It follows from Proposition~\ref{pr:fk} that the entry $t_{ij}$ of the
transition matrix, where $i=0, 1, \ldots, p-1$ and $j=0, 1, \ldots,
p-1$ is equal to
\begin{multline*}
t_{ij}=(-1)^j\binom{i+j}{j}=(-1)^j\frac{(i+j)(i+j-1)\cdots(i+j-j+1)}{j!}=\\
\frac{(-i-j)(-i-j+1)\cdots(-i-1)}{j!}=\\
\frac{(p-i-j)(p-i-j+1)\cdots(p-i-1)}{j!}=\binom{p-1-i}{j},\end{multline*}
which proves the first claim of the proposition.

In order to prove the second claim, we have to show that the
product
\[\left(\binom{p-1-i}{j}\right)_{i,j=0}^{p-1}\cdot\left(\binom{j}{p-1-i}\right)_{i,j=0}^{p-1}\]
is equal to the identity matrix. The general entry of the product is
equal to
\[a_{ij}=\sum_{k=0}^{p-1}\binom{p-1-i}{k}\cdot\binom{j}{p-1-k}=
\binom{p-1-i+j}{p-1}.\]
But $\binom{x+p-1}{p-1}=\frac{(x+p-1)(x+p-2)\cdots(x+1)}{(p-1)!}$ is
equal to one for $x=0$ and is equal to zero for $x\ne 0$, which shows
that the product is equal to the identity matrix.
\end{proof}

\begin{examp}
\label{ex:admachmatrix}
In the case $p=2$, the transition matrix is
$T=\left(\begin{array}{cc}1 & 1\\ 1 & 0\end{array}\right)$, and its
inverse is $T^{-1}=\left(\begin{array}{cc} 0 & 1\\ 1 &
    1\end{array}\right)$. Let us use this to find the matrix
recursions for some self-similar groups acting on the binary tree in
the new basis $\mathsf{B}=\{b_0, b_1\}$.

For the adding machine (see Examples~\ref{ex:admach},
\ref{ex:admach2}, \ref{ex:admach3}), we have
\[\stencil_2(a)=\left(\begin{array}{cc} 0 & 1\\ 1 &
    1\end{array}\right)\cdot
\left(\begin{array}{cc}0 & a\\ 1 &
    0\end{array}\right)\cdot\left(\begin{array}{cc} 1 & 1\\ 1 &
    0\end{array}\right)=
\left(\begin{array}{cc}1 & 1\\ 1+a & 1\end{array}\right).\]
\end{examp}

\begin{proposition}
The matrix of the binary adding machine $\pi_\infty(a)$ in the basis $\mathsf{B}_\infty$
is the infinite Jordan cell
\[\left(\begin{array}{ccccc} 1 & 1 & 0 & 0 & \ldots\\
0 & 1 & 1 & 0 & \ldots\\ 0 & 0 & 1 & 1 & \ldots\\ 0 & 0 & 0 & 1 &
\ldots\\
\vdots & \vdots & \vdots & \vdots & \ddots\end{array}\right).\]
\end{proposition}

\begin{proof}
Let us prove the statement by induction, using the matrix recursion
from Example~\ref{ex:admachmatrix}.
The matrix of $a$ on $V_0$ is $(1)$. The four 2-decimations
$\stencil_2(J)_{ij}$ of the
Jordan cell are
\[\stencil_2(J)_{00}:\left(\begin{array}{ccccc} \raisebox{.5pt}{\textcircled{\raisebox{-.9pt} {1}}} & 1 & \raisebox{.5pt}{\textcircled{\raisebox{-.9pt} {0}}} & 0 & \ldots\\
0 & 1 & 1 & 0 & \ldots\\ \raisebox{.5pt}{\textcircled{\raisebox{-.9pt} {0}}} & 0 & \raisebox{.5pt}{\textcircled{\raisebox{-.9pt} {1}}} & 1 & \ldots\\ 0 & 0 & 0 & 1 &
\ldots\\
\vdots & \vdots & \vdots & \vdots & \ddots\end{array}\right)\mapsto
\left(\begin{array}{cccc} 1 & 0 & 0 &  \ldots\\
0 & 1 & 0 & \ldots\\
0 & 0 & 1 & \ldots\\
\vdots & \vdots & \vdots & \ddots\end{array}\right),\]
\[\stencil_2(J)_{01}:\left(\begin{array}{ccccc} 1 &
    \raisebox{.5pt}{\textcircled{\raisebox{-.9pt} {1}}} & 0 &\raisebox{.5pt}{\textcircled{\raisebox{-.9pt} {0}}} & \ldots\\
0 & 1 & 1 & 0 & \ldots\\ 0 &
\raisebox{.5pt}{\textcircled{\raisebox{-.9pt} {0}}} & 1 & \raisebox{.5pt}{\textcircled{\raisebox{-.9pt} {1}}} & \ldots\\ 0 & 0 & 0 & 1 &
\ldots\\
\vdots & \vdots & \vdots & \vdots & \ddots\end{array}\right)\mapsto
\left(\begin{array}{cccc} 1 & 0 & 0 &  \ldots\\
0 & 1 & 0 & \ldots\\
0 & 0 & 1 & \ldots\\
\vdots & \vdots & \vdots & \ddots\end{array}\right),\]
\[\stencil_2(J)_{10}:
\left(\begin{array}{ccccc}
1 & 1 & 0 & 0 & \ldots\\
\raisebox{.5pt}{\textcircled{\raisebox{-.9pt} {0}}} & 1 &
\raisebox{.5pt}{\textcircled{\raisebox{-.9pt} {1}}} & 0 & \ldots\\
0 & 0 & 1 & 1 & \ldots\\ \raisebox{.5pt}{\textcircled{\raisebox{-.9pt} {0}}} & 0 & \raisebox{.5pt}{\textcircled{\raisebox{-.9pt} {0}}} & 1 & \ldots\\
\vdots & \vdots & \vdots & \vdots & \ddots\end{array}\right)\mapsto
\left(\begin{array}{cccc} 0 & 1 & 0 & \ldots\\
0 & 0 & 1 & \ldots\\
0 & 0 & 0 & \ldots\\
\vdots & \vdots & \vdots & \ddots\end{array}\right),\]
\[\stencil_2(J)_{11}:
\left(\begin{array}{ccccc}
1 & 1 & 0 & 0 & \ldots\\
0 & \raisebox{.5pt}{\textcircled{\raisebox{-.9pt} {1}}} & 1 &
\raisebox{.5pt}{\textcircled{\raisebox{-.9pt} {0}}} & \ldots\\
0 & 0 & 1 & 1 & \ldots\\ 0 & \raisebox{.5pt}{\textcircled{\raisebox{-.9pt} {0}}} & 0 & \raisebox{.5pt}{\textcircled{\raisebox{-.9pt} {1}}}  & \ldots\\
\vdots & \vdots & \vdots & \vdots & \ddots\end{array}\right)\mapsto
\left(\begin{array}{cccc} 1 & 0 & 0 &  \ldots\\
0 & 1 & 0 & \ldots\\
0 & 0 & 1 & \ldots\\
\vdots & \vdots & \vdots & \ddots\end{array}\right),\]
which agrees with the recursion $\stencil_2(a)=\left(\begin{array}{cc}1 & 1
    \\ 1+a & 1\end{array}\right)$.
\end{proof}

\begin{lemma}
\label{lem:5.7}
If $g\in\autxs$ satisfies the wreath recursion $g=(g_0, g_1)$, then
its matrix recursion in the basis $\mathsf{B}_\infty$ (over the field
$\fieldp$ for $p=2$) is
\[\stencil_2(g)=\left(\begin{array}{cc} g_1 & 0 \\ g_0+g_1 &
    g_0\end{array}\right).\]
If it satisfies $g=\sigma(g_0, g_1)$, then
\[\stencil_2(g)=\left(\begin{array}{cc} g_0 & g_0\\ g_0+g_1 &
    g_0\end{array}\right).\]
\end{lemma}

\begin{proof}
We have, in the first case,
\[\stencil_2(g)=\left(\begin{array}{cc}0 & 1\\ 1 & 1\end{array}\right)
\left(\begin{array}{cc} g_0 & 0 \\ 0 &
    g_1\end{array}\right)\left(\begin{array}{cc}
1 & 1 \\ 1  & 0\end{array}\right)=\left(\begin{array}{cc} g_1 & 0\\
g_0+g_1 & g_0\end{array}\right).\]
In the second case:
\[\stencil_2(g)=\left(\begin{array}{cc}0 & 1\\ 1 & 1\end{array}\right)
\left(\begin{array}{cc} 0 & g_1 \\ g_0 &
    0\end{array}\right)\left(\begin{array}{cc}
1 & 1 \\ 1  & 0\end{array}\right)=\left(\begin{array}{cc} g_0 & g_0\\
g_0+g_1 & g_0\end{array}\right).\]
\end{proof}

\begin{examp}
\label{ex:grimatrices}
It follows from Lemma~\ref{lem:5.7} that
the matrix recursion for the generators of the Grigorchuk group (see
Example~\ref{ex:grigorchuk}) in the basis $\mathsf{B}$ is
\[\stencil_2(a)=\left(\begin{array}{cc}1 & 1\\ 0 &
    1\end{array}\right),\quad
\stencil_2(b)=\left(\begin{array}{cc}c & 0\\ a+c &
    a\end{array}\right),\]
\[\stencil_2(c)=\left(\begin{array}{cc}d & 0\\ a+d &
    a\end{array}\right),\quad
\stencil_2(d)=\left(\begin{array}{cc}b & 0\\ 1+b &
    1\end{array}\right).\]
See a visualization of the matrices $b, c, d$ on Figure~\ref{fig:bcd},
where black pixels correspond to ones, and white pixels to zeros.

\begin{figure}
\centering
\includegraphics[width=4in]{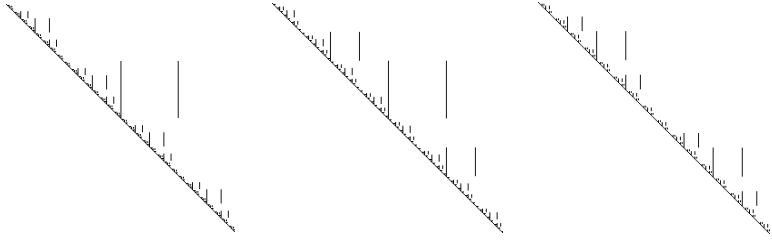}
\caption{The matrices of the generators $b, c, d$ of the Grigorchuk group}
\label{fig:bcd}
\end{figure}
\end{examp}

Denote $U_n=\langle b_0, b_1, \ldots, b_n\rangle=\langle e_0, e_1,
\ldots, e_n\rangle<C(X^\omega, \fieldp)$.

\begin{proposition}
Each space $U_n$ is $\sylow$-invariant,
and the kernel in $\sylow$ of the restriction of $\pi_\infty$ onto
$U_{p^{n-1}}=\langle
b_0, b_1, \ldots, b_{p^{n-1}}\rangle$ coincides with the kernel of
$\pi_n$. In other words, restriction of $\pi_\infty$ onto
$U_{p^{n-1}}$ defines a faithful representation of $\sylown$.
\end{proposition}

\begin{proof}
The subspace $\langle b_0, b_1, \ldots, b_{p^{n-1}}\rangle<
C(X^\omega, \fieldp)$ is equal to the span of the product $V_{n-1}\cdot
e_{p^{n-1}}$, where $e_{p^{n-1}}$ is the function on $X^\omega$ given
by
\[e_{p^{n-1}}(x_1x_2\ldots)=x_n,\]
according to~\eqref{eq:en}. In other words, it is the tensor product
$V_{n-1}\otimes\langle e_1\rangle$, where $e_1(x)=x$.

Suppose that $g\in G$ belongs to the kernel of the restriction of
$\pi_\infty$ onto $V_{n-1}\otimes\langle e_1\rangle$. Then for every
$v\in X^{n-1}$ we have $\pi_1(g)(\delta_v)=\delta_v$, since
$\delta_v\in V_{n-1}$. Then
\[\pi_\infty(g)(\delta_v\otimes e_1)=\delta_v\otimes (e_1\circ
\pi_1(g|_v)^{-1}),\]
hence $\pi_1(g|_v)$ is identical for every $v\in X^{n-1}$. It follows
that $g$ acts trivially on $X^n$, i.e., that $\pi_n(g)$ is trivial.
\end{proof}

Thus, we get a faithful representation of $\sylown=\wr_{k=1}^nC_p$ by uni-triangular matrices of
dimension $p^{n-1}+1$. Note that this is the smallest possible
dimension for a faithful representation, since the nilpotency class
of $\sylown$ is equal to $p^{n-1}$, while the nilpotency
class of the group of uni-triangular matrices of dimension $d$ is
equal to $d-1$.

\subsection{The first diagonal}

Let $\alpha:\sylow\arr\fieldp^\omega$ be the abelianization homomorphism given
by~\eqref{eq:alphamap}. We write \[\alpha(g)=(\alpha_0(g), \alpha_1(g), \ldots).\]

If $A=(a_{ij})_{i, j=0}^\infty$ is an infinite matrix, then its
\emph{first diagonal} is the sequence $(a_{01}, a_{12}, a_{23},
\ldots)$, i.e., the first diagonal above the main diagonal of $A$.

\begin{theorem}
\label{th:diagonal}
Let $g\in G$, and let $A_g=(a_{ij})_{i, j=0}^\infty$ be the matrix of
$\pi_\infty(g)$ in the basis $\mathsf{B}_\infty$, constructed in the
previous section. Let $(s_1, s_2,
\ldots)=(a_{01}, a_{12}, \ldots)$ be the first diagonal of $A_g$. Then
\[s_n=\alpha_k(g),\]
where $p^k$ is the maximal power of $p$ dividing $n$.
\end{theorem}

For example, if $p=2$, and $\alpha(g)=(a_0, a_1, a_2, \ldots)$, then
the first diagonal of $A_g$ is \[a_0, a_1, a_0, a_2, a_0, a_1, a_0,
a_3, a_0, a_1, a_0, a_2, \ldots\]

\begin{proof}
The first diagonal of a product of two upper uni-triangular matrices
$A$ and $B$ is equal to the sum of the first diagonals of the matrices
$A$ and $B$. It follows that it is enough to prove the theorem for
rooted automorphisms of $\autxs$ (i.e., automorphisms $g$ such that
$g|_v$ is trivial for all non-empty words $v\in X^*$) and for automorphisms acting
trivially on the first level $X$.

If the automorphism is rooted, then it is a power of the automorphism
\[a:x_1x_2\ldots x_n\mapsto (x_1+1)x_2\ldots x_n.\]

It follows from the definition of the basis $\mathsf{B}$ that the matrix of
$\pi_\infty(a)$ is the block-diagonal matrix consisting of the Jordan
cells of size $p$. Consequently, its first diagonal is the periodic
sequence of period $(1, 1, \ldots, 1, 0)$ of length $p$. Hence, the
first diagonal of $a^s$ is $(s, s, \ldots, s, 0)$ repeated
periodically. This proves the statement of the theorem for the
automorphisms of the form $a^s$.

Suppose that $g$ acts trivially on the first level of the tree. Then
its matrix recursion in the basis $\{\delta_x\}_{x\in X}$ is the
diagonal matrix with the entries $g|_x$ on the diagonal. It follows
that the matrix recursion for $g$ in the basis $\{b_i\}_{i=0}^{p-1}$ is
equal to the product of the matrices
\begin{multline*}
\left(\begin{array}{ccccc}0 & 0 & 0 & \ldots & \binom{p-1}{p-1}\\
0 & 0 & 0 & \ldots & \binom{p-1}{p-2}\\
\vdots & \vdots & \vdots & \ddots & \vdots \\
0 & 0 & 1 & \ldots & \binom{p-1}{2}\\
0 & 1 & 2 & \ldots & \binom{p-1}{1}\\
1 & 1 & 1 & \ldots & 1\end{array}\right)\cdot
\left(\begin{array}{ccccc}g|_0 & 0 & 0 & \ldots  & 0\\
0 & g|_1 & 0 & \ldots & 0\\
0 & 0 & g|_2 & \ldots & 0\\
\vdots & \vdots & \vdots & \ddots & \vdots \\
0 & 0 & 0 & \ldots & g|_{p-1}\end{array}\right)\cdot\\
\left(\begin{array}{ccccc}1 & \binom{p-1}{1} & \binom{p-1}{2} & \ldots
    & \binom{p-1}{p-1}\\
1 & \binom{p-2}{1} & \binom{p-2}{2} & \ldots & 0\\
\vdots & \vdots & \vdots & \ddots & \vdots \\
1 & 2 & 1 & \ldots & 0\\
1 & 1 & 0 & \ldots & 0\\
1 & 0 & 0 & \ldots & 0\end{array}\right).
\end{multline*}
It is easy to see that the entries on the first diagonal above the
main diagonal of the product are equal to zero, and that the entry in
the left bottom corner is equal to $g|_0+g|_1+\cdots+g|_{p-1}$.

If we apply the stencil map
\[\stencil_p(s_1, s_2, \ldots)=((s_1, s_{1+p}, \ldots), (s_2, s_{2+p},
\ldots), \ldots, (s_p, s_{2p}, \ldots))=(w_1, w_2, \ldots, w_p)\]
to the first diagonal of $A_g$, then $w_1, w_2, \ldots, w_{p-1}$ are the
main diagonals of the $p$-decimations $B_1, B_2, \ldots, B_{p-1}$ of
$A_g$, where $(B_1, B_2, \ldots, B_{p-1})$ is the first diagonal above
the main diagonal of $\stencil_p(A_g)$. The sequence $w_p$ is the
first diagonal of the entry in the lower left corner of
$\stencil_p(A_g)$. It follows that $(s_1, s_2, \ldots)$ is of the form
$(0, 0, \ldots, s_1', 0, 0, \ldots, s_2', \ldots)$, where there are
$p-1$ zeros at the beginning and between the entries $s_i'$, and
$(s_1', s_2', \ldots)$ is the first diagonal of
$\pi_\infty(g|_0+g|_1+\cdots+g|_{p-1})$ in the basis $(b_0, b_1,
\ldots, b_{p-1})$. This provides us with an inductive proof of the
statement of the theorem.
\end{proof}

\begin{examp}
Consider the matrices $a, b, c, d$ generating the Grigorchuk group, as
it is described in Example~\ref{ex:grimatrices}. It follows directly
from the description of the action of the elements $a, b, c, d$ on the
rooted tree (see Example~\ref{ex:grigorchuk}) that
\begin{eqnarray*}
\alpha(a) &=& (1, 0, 0, 0, 0, 0, 0, 0, 0, 0, \ldots),\\
\alpha(b) &=& (0, 1, 1, 0, 1, 1, 0, 1, 1, 0, \ldots),\\
\alpha(c) &=& (0, 1, 0, 1, 1, 0, 1, 1, 0, 1, \ldots),\\
\alpha(d) &=& (0, 0, 1, 1, 0, 1, 1, 0, 1, 1, \ldots).
\end{eqnarray*}
It follows from Theorem~\ref{th:diagonal} that the first diagonal of
$a$ is $(1, 0, 1, 0, \ldots)$.
The first diagonal of $b$ is $(s_1, s_2, \ldots)$ where
\[s_{2^{3k}(2m+1)}=0, \quad s_{2^{3k+1}(2m+1)}=1, \quad s_{2^{3k+2}(2m+1)}=1,\]
where $k, m\ge 0$ are integers.
\end{examp}

Sequence $s_n$ from the previous example is a \emph{Toeplitz
sequence}, see~\cite{toeplitz,prodingerurbanek}
and~\cite[Exercise~10.11.42]{allouchshallit:sequences}. If $G$ is a finitely
generated finite state self-similar group, then the sequence $\alpha_k(g)$ in
Theorem~\ref{th:diagonal} is eventually periodic for every $g\in G$,
and then the sequence $s_n$ is Toeplitz.

\subsection{Generating series}

Let $w=(a_0, a_1, \ldots)\in\Bbbk^\omega$, and consider the
corresponding formal power series
$a_0+a_1x+a_2x^2+\cdots=G_w(x)\in\Bbbk[[x]]$. It is easy to see that
if
\[\stencil_d(w)=(w_0, w_1, \ldots, w_{d-1}),\]
then
\[G_w(x)=G_{w_0}(x^d)+xG_{w_1}(x^d)+\cdots+x^{d-1}G_{w_{d-1}}(x^d).\]

Note that if $\Bbbk=\mathbb{F}_p$ and $d=p$, then we get
\[G_w(x)=(G_{w_0}(x))^p+x(G_{w_1}(x))^p+\cdots+x^{p-1}(G_{w_{p-1}}(x))^p.\]

We have the following characterization of automatic sequences, due to
Christol, see~\cite[Theorem~12.2.5]{allouchshallit:sequences}.

\begin{theorem}
Let $\Bbbk$ be a finite field of characteristic $p$.
Then a sequence $w\in\Bbbk^\omega$ is $p$-automatic if
and only if the generating series $G_w(x)$ is algebraic over $\Bbbk(x)$.
\end{theorem}

Similarly, if $A=(a_{ij})_{i, j=0}^\infty$ is a matrix over a field
$\Bbbk$, then we can consider the formal series $G_A(x, y)=\sum_{i, j=0}^\infty
a_{ij}x^iy^j\in\Bbbk[[x,
y]]$. If \[\stencil_d(A)=\left(A_{ij}\right)_{i, j=0}^{d-1},\]
then
\[G_A=\sum_{i, j=0}^{d-1}x^iy^jG_{A_{ij}}(x^d, y^d).\]

We also have a complete analog of Christol's theorem for matrices,
see~\cite[Theorem~14.4.2]{allouchshallit:sequences}.

\begin{theorem}
Let $\Bbbk$ be a finite field of characteristic $p$. Then a matrix
$A=(a_{ij})_{i, j=0}^\infty$ is $p$-automatic if and only if the
series $G_A$ is algebraic over $\Bbbk(x, y)$.
\end{theorem}

In the case when $A$ is triangular, it may be natural to use the
generating function
\[T_A(t, s)=\sum a_{ij}t^{j-i}s^i,\]
so that $T_A(t, s)=H_0(s)+H_1(s)t+H_2(s)t^2+\cdots,$
where $H_i(y)$ are generating functions of the diagonals of $A$. Note
that $T_A$ and $G_A$ are related by the formula:
\[T_A(t, s)=G_A(s/t, t)\]

\begin{examp}
Consider the generators of the Grigorchuk group $a, b, c, d$ given by
the matrices from Example~\ref{ex:grimatrices}. Let $A(x, y), B(x, y),
C(x, y)$, and $D(x, y)$ be the corresponding generating series. Note
that the generating series of the unit matrix is $I(x,
y)=\frac{1}{1+xy}=\frac{1}{1+s}$.

It follows from the recursions in Example~\ref{ex:grimatrices} that
\[A=I(x^2, y^2)+yI(x^2, y^2)+xyI(x^2,
y^2)=\frac{1+y+xy}{1+x^2y^2}=\frac{1}{1+s}+\frac{t}{1+s^2}.\]

\[B=C^2+x(A^2+C^2)+xyA^2=
(1+x)C^2+(x+xy)A^2
\]

\[C=(1+x)D^2+(x+xy)A^2
\]

\[D=(1+x)B^2+(x+xy)I^2.\]

Let us make a substitution $A=y\tilde A+\frac{1}{1+xy}$,
$B=y\tilde B+\frac{1}{1+xy}$, $C=y\tilde C+\frac{1}{1+xy}$, and
$D=y\tilde D+\frac{1}{1+xy}$. Note that $I=\frac{1}{1+xy}$, so we set
$\tilde I=0$. The series $\tilde A, \tilde B, \tilde C, \tilde D$ are
the generating series of the matrices obtained from the matrices $a,
b, c, d$ by removing the main diagonal, and shifting all columns to
the left by one position.

We have then
\[\tilde A=\frac{1}{1+x^2y^2}=
\frac{1}{1+s^2}\]
and
\[y\tilde B+\frac{1}{1+xy}=(1+x)\left(y^2{\tilde C}^2+\frac{1}{1+x^2y^2}\right)+(x+xy)\left(y^2{\tilde A}^2+\frac{1}{1+x^2y^2}\right),\]
hence
\[y\tilde B=(1+x)y^2{\tilde C}^2+(x+xy)y^2{\tilde A}^2,\]
hence
\[\tilde B=(y+xy){\tilde C}^2+\frac{xy+xy^2}{1+x^4y^4}=(s+t){\tilde C}^2+\frac{s+ts}{1+s^4}.\]
Similarly,
\[\tilde C=(y+xy){\tilde D}^2+\frac{xy+xy^2}{1+x^4y^4}=(s+t){\tilde D}^2+\frac{s+ts}{1+s^4},\]
and
\[\tilde D=(y+xy){\tilde B}^2=(s+t){\tilde B}^2.\]
Let us denote $F=\frac{s+ts}{1+s^4}$.
Then $\tilde B$, $\tilde C$, and $\tilde D$ are solutions of the equations
\[(s+t)^7{\tilde B}^8+\tilde B+(s+t)F^2+F=0\]
\[(s+t)^7{\tilde C}^8+\tilde C+(s+t)^3F^4+F=0,\]
\[(s+t)^7{\tilde D}^8+\tilde D+(s+t)^3F^4+(s+t)F^2=0.\]

Substituting $t=0$, we get equations for the generating functions
$B_1(s), C_1(s), D_1(s)$ of the first diagonals above the
main in the matrices $b, c, d$:
\begin{equation}
\label{eq:B1}
s^7B_1^8+B_1+\frac{s^3}{1+s^8}+\frac{s}{1+s^4}=0,
\end{equation}
\begin{equation}
\label{eq:C1}
s^7C_1^8+C_0+\frac{s^7}{1+s^{16}}+\frac{s}{1+s^4}=0,
\end{equation}
\begin{equation}
\label{eq:D1}
s^7D_1^8+D_1+\frac{s^7}{1+s^{16}}+\frac{s^3}{1+s^8}=0.
\end{equation}
\end{examp}

Denote
\[J=\left(\begin{array}{ccccc}
0 & 1 & 0 & 0 & \ldots\\
0 & 0 & 1 & 0 & \ldots\\
0 & 0 & 0 & 1 & \ldots\\
\vdots & \vdots & \vdots & \vdots & \ddots\end{array}\right)\]
Then every upper uni-triangular matrix $M$ can be written as
\begin{equation}
\label{eq:JM}
I+D_1(M)J+D_2(M)J^2+D_3(M)J^3+\cdots,
\end{equation}
where $D_i(M)$ are diagonal matrices whose main diagonals are equal to the
$i$th diagonals of $M$.

The generating series $T_M(t, s)$ is equal to
\begin{equation}
\label{eq:poserst}
\frac{1}{1+s}+M_1(s)t+M_2(s)t^2+M_3(s)t^3+\cdots,
\end{equation}
where $M_i(s)$ is the usual generating series of the main diagonal of
$D_i(M)$. Addition and multiplication of diagonal matrices $\diag(a_0,
a_1, \ldots)$ corresponds to the usual addition and the
\emph{Hadamard} (coefficient-wise) multiplication of the power series
$a_0+a_1s+a_2s^2+\cdots$. Note that we have
\[J\diag(a_0, a_1, \ldots)=\diag(a_1, a_2, \ldots)J,\]
which gives an algebraic rule for multiplication of the power
series~\eqref{eq:poserst} corresponding to multiplication of matrices.

Namely, we can replace the matrix $M$ by the formal power
series~\eqref{eq:poserst}, where the series in the variable $s$ are
added and multiplied coordinate-wise, while the series in the variable
$t$ are multiplied in the usual (though non-commutative) way subject
to the relation
\[t(a_0s^0+a_1s+a_2s^2+\cdots)=(a_1s^0+a_2s+a_3s^2+\cdots)t.\]

Let $A=\left(a_{ij}\right)_{i, j=0}^\infty$ be a matrix, and denote by
$\Delta_i(M)=(a_{0, i}, a_{1, i+1}, a_{2, i+2}, \ldots)$ the
sequence equal to the $i$th diagonal of $A$.

Let $d>2$, $k\in\{0, 1, \ldots, d-1\}$, $n\ge 1$ be integers, and let
$k+n=dq+r$ for $q\in\mathbb{Z}$ and $r\in\{0, 1, \ldots, d-1\}$. Then
the sequence
$\stencil_d(\Delta_n(M))_k=(a_{k, dq+r}, a_{k+d, d(q+1)+r}, a_{k+2d,
  d(q+2)+r}, \ldots)$ is equal to the $q$th diagonal of the matrix
$(\stencil_d{M})_{k, r}$.

\begin{examp}
Consider, as in the previous examples, the matrices $a, b, c, d$
generating the Grigorchuk group. Denote $A_n=\Delta_n(a), B_n=\Delta_n(b),
C_n=\Delta_n(c), D_n=\Delta_n(d)$. Note that $A_0=(1, 1, \ldots)$,
$A_1=(1, 0, 1, 0, \ldots)$, and $A_n=0$ for all $n\ge 2$. We have
\[\stencil_2(B_{2n})=\left(C_n, A_n\right)\]
\[\stencil_2(B_{2n+1})=\left(0, C_{n+1}+A_{n+1}\right).\]
Similarly,
\[\stencil_2(C_{2n})=\left(D_n, A_n\right)\]
\[\stencil_2(C_{2n+1})=\left(0, A_{n+1}+D_{n+1}\right),\]
and
\[\stencil_2(D_{2n})=\left(B_n, I_n\right),\]
\[\stencil_2(D_{2n+1})=\left(0, I_{n+1}+B_{n+1}\right).\]

These stencil recursions give us recursive formulas for the
corresponding generating functions, which we will denote $B_n(s)$,
$C_n(s)$, $D_n(s)$. Recall that $A_n(s)=0$ for $n\ge 2$, $I_n(s)=0$ for
$n\ge 1$, $A_1(s)=\frac{1}{1+s^2}$, and $A_0(s)=I_0(s)=\frac{1}{1+s}$.
\begin{eqnarray*}
B_{2n}(s) &=& C_n^2+sA_n^2\\
B_{2n+1}(s) &=& s(C_{n+1}^2+A_{n+1}^2)\\
C_{2n}(s) &=& D_n^2+sA_n^2\\
C_{2n+1}(s) &=& s(D_{n+1}^2+A_{n+1}^2)\\
D_{2n}(s) &=& B_n^2+sI_n^2\\
D_{2n+1}(s) &=& s(I_{n+1}^2+B_{n+1}^2)
\end{eqnarray*}

Note that iterations of the map $2n\mapsto n, 2n+1\mapsto n+1$ on the
set of non-negative integers are
attracted to two fixed points $0\mapsto 0$ and $1\mapsto 1$. Consequently, we
get the following
\begin{proposition}
For every $n\ge 1$ the generating functions $B_n(s), C_n(s), D_n(s)$
are of the
form \[\frac{p_0(s)}{1+s^{2^k}}+p_1(s)B_1(s)^{2^l},\]
where $k, l\ge 0$ are integers, and $p_0(s), p_1(s)$ are polynomials
over $\mathbb{F}_2$.
\end{proposition}
\end{examp}

\subsection{Principal columns}
\label{ss:principal}

Let $g\in\sylow$, let $A_g$ be the matrix of $\pi_\infty(g)$ in the
basis $\mathsf{E}_\infty$ of monomials. Recall that we number the columns
and rows of the matrix $A_g$ starting from zero.

\begin{proposition}
\label{pr:principal}
Every entry $a_{i, j}$ of the matrix $A_g$ is a polynomial function
(not depending on $g$) of
the entries of the columns number $p, p^2, \ldots,
p^{\lfloor\log_p j\rfloor}$.

The same statement is true for the matrices of $\pi_\infty(g)$ in the
basis $\mathsf{B}_\infty$.
\end{proposition}

\begin{proof}
Let
\[[a_0, a_1(x_1), a_2(x_1, x_2), \ldots]\] be the tableau of $g$, as
in Subsection~\ref{ss:sylow}.

Recall that $e_{p^n}$ is the monomial $x_n$. Consequently,
\[\pi_\infty(g)(e_{p^n})=x_n+u_n(x_1, x_2, \ldots, x_{n-1}).\]
It follows that the entries of column number $p^n$ of the matrix $A_g$
above the main diagonal are the coefficients of the representation of $u_n$ as a linear
combination of monomials $e_i$ for $i<p^n$, i.e., are the coefficients
of the polynomial $u_n$. (The entries below the diagonal are zeros,
and the entry on the diagonal is equal to one, of course.)

If $k$ is a natural number such that $k<p^{n+1}$, then
$e_k=x_1^{r_1}x_2^{r_2}\cdots x_n^{r_n}$, where $r_n, r_{n-1}, \ldots,
r_1$ are the digits of the representation of $k$ in the base $p$
numeration system, and
\[\pi_\infty(g)(e_k)=(x_1+u_1)^{r_1}(x_2+u_2(x_1))^{r_2}\cdots(x_n+u_n(x_1,
x_2, \ldots, x_{n-1})^{r_n},\]
which implies the statement of the proposition.

For the basis $\mathsf{B}_\infty$ we have $b_{p^n}=-\binom{x_n+1}{1}=-x_n-1$,
and a similar proof works.
\end{proof}

\begin{defi}
Columns number $p^n$, $n=0, 1, 2, \ldots$, of the matrix $A_g$ are called the
\emph{principal columns} of $A_g$.
\end{defi}

\begin{examp}
Let, for $p=2$, the first four principal columns of the matrix $A_g$,
$g\in\sylow$ be $(a_{01}, 1, 0, \ldots)^\top$, $(a_{02}, a_{12}, 1, 0,
\ldots)^\top$, and $(a_{04}, a_{14}, a_{24}, a_{34}, 1, 0,
\ldots)^\top$. Then the columns number 3, 5, 6, and 7 (when numeration
of the columns starts from zero) are
\[\left(\begin{array}{c}
a_{01}a_{02}\\
a_{01}a_{12}+a_{12}+a_{02}\\
a_{01}\\
1
\end{array}\right),
\left(\begin{array}{c}
a_{01}a_{04}\\
a_{01}a_{14}+a_{14}+a_{04}\\
a_{01}a_{24}\\
a_{24}+a_{34}a_{01}+a_{34}\\
a_{01}\\
1
\end{array}\right),\]
\[\left(
\begin{array}{c}
a_{02}a_{04}\\
a_{02}a_{14}+a_{12}a_{04}+a_{12}a_{14}\\
a_{04}+a_{24}a_{02}+a_{24} \\
a_{14}+a_{34}+a_{02}a_{34}+a_{12}a_{34}+a_{12}a_{24} \\
a_{02} \\
a_{12} \\
1
\end{array}\right),
\]
and
\[
\left(\begin{array}{c}
a_{01}a_{02}a_{04}\\
(a_{04}+a_{14}a_{01}+a_{14})(a_{02}+a_{12})+a_{01}a_{12}a_{04}\\
a_{01}a_{04}+a_{01}a_{24}a_{02}+a_{01}a_{24} \\
(a_{02}+a_{12}+1)(a_{24}+a_{34}+a_{01}a_{34})+
a_{01}(a_{14}+a_{12}a_{24})+a_{04}+a_{14}
\\
 a_{01}a_{02}\\
a_{02}+a_{12}+a_{01}a_{12} \\
a_{01}\\
1
\end{array}\right),
\]
respectively.
\end{examp}

\subsection{Uniseriality}

Recall that the basis $\mathsf{E}_\infty$ of $C(X^\omega, \fieldp)$
consists of monomial functions $e_0, e_1, \ldots$, where
\[e_{k_0+k_1p+k_2p^2+\cdots}(x_1, x_2,
\ldots)=x_1^{k_0}x_2^{k_1}x_3^{k_2}\cdots,\]
where $k_i\in\{0, 1, \ldots, p-1\}$ are almost all equal to zero.
Let us call, following~\cite{kalouj:sylowfin} (though we use a
slightly different definition),
$n=k_0+k_1p+k_2p^2+\cdots$ the \emph{height} of the monomial
$e_n=x_1^{k_0}x_2^{k_1}x_3^{k_2}\cdots$.

Height of a reduced polynomial $f$ is defined as the maximal height of its
monomials, and is denoted $\gamma(f)$. We define $\gamma(0)=-1$ (note that our
definition is different from the definition of Kaloujnine, which uses
$\gamma(f)+1$, so that height of $0$ is zero).

Let us describe, following~\cite{tullioleonovscarabottitolli}, an algorithm for computing height of a
function $f\in\fieldp^{X^n}$. Let $\mathsf{B}$ be the basis
$b_0<b_1<\ldots<b_{p-1}$ of $\fieldp^X$, constructed
in~\ref{ss:polynombases}.

Let $\{\delta_0', \delta_1', \ldots, \delta_{p-1}'\}$ and
$\{b_0', b_1', \ldots, b_{p-1}'\}$ be the bases of the space of functionals
$\left(\fieldp^X\right)'$ dual to the bases $\{\delta_0, \delta_1,
\ldots, \delta_{p-1}\}$ and $\mathsf{B}$, respectively, i.e.,
$\delta_i'$ and $b_i'$ are defined by the condition
\[\langle\delta_i'|\delta_j\rangle=\langle b_i'|b_i\rangle=\delta_{i,
  j}\]
for all $0\le i, j\le p-1$.

Then we have
\[\langle\delta_x'|f\rangle=f(x),\]
for all $f\in\fieldp^X$ and $x\in X=\fieldp$.

It follows from Proposition~\ref{pr:transitionmatrixb} and elementary
linear algebra that the transition matrix from the basis $\{\delta_0',
\delta_1', \ldots, \delta_{p-1}'\}$ to $\{b_0', b_1', \ldots,
b_{p-1}'\}$ is the matrix transposed to the matrix $T^{-1}$ of
Proposition~\ref{pr:transitionmatrixb}, i.e., the matrix
\[\left(\begin{array}{cccccc} 0 & \ldots & 0 & 0 & 0 & 1\\
0 & \ldots & 0 & 0 & 1 & 1\\
0 & \ldots & 0 & 1 & 2 & 1\\
0 & \ldots & 1 & 3 & 3 & 1\\
\vdots & \ldots & \vdots & \vdots & \vdots & \vdots\\
\binom{p-1}{p-1} & \ldots & \binom{p-1}{3} & \binom{p-1}{2} &
\binom{p-1}{1} & \binom{p-1}{0}
\end{array}\right).\]
In other words,
\begin{equation}\label{eq:bkprime}
b_k'=\sum_{l=0}^{p-1}\binom{l}{p-1-k}\delta_l',
\end{equation}
so that
\[\langle b_k'|f\rangle=\sum_{l=0}^{p-1}\binom{l}{p-1-k}f(l).\]
For instance, $\langle b_{p-1}'|f\rangle=\sum_{x\in X}f(x)$.

Define linear maps $R_k:\fieldp^{X^{n+1}}\arr\fieldp^{X^n}$, $k=0, 1,
\ldots, p-1$, as the linear extension of the map
\[b_{i_1}\otimes \cdots b_{i_n}\otimes b_{i_{n+1}}\mapsto
\langle b_k'|b_{i_{n+1}}\rangle\cdot b_{i_1}\otimes \cdots b_{i_n}.\]

In other terms, the map $R_k$ is given by
\[R_k(f)(x_1, x_2, \ldots, x_n)=\langle b_k'|f(x_1, x_2, \ldots, x_n,
x)\rangle,\]
where $f(x_1, x_2, \ldots, x_n, x)\in\fieldp^{X^{n+1}}$ on the
right-hand side of the equality is treated as a function of $x$ for
every choice of  $(x_1, x_2, \ldots, x_n)\in\fieldp^{X^n}$.

Using~\eqref{eq:bkprime}, we see
that $R_k$ can be computed using the formula
\[R_k(f)(x_1, x_2, \ldots, x_n)=\sum_{x=0}^{p-1}\binom{x}{p-1-k}f(x_1,
x_2, \ldots, x_n, x).\]

\begin{proposition}
\label{pr:heightp}
Let $f\in\fieldp^{X^n}$. Define $j_n$ as the maximal value of $j=0, 1,
\ldots, p-1$ such that $R_j(f)\ne 0$, and then define inductively $j_k$ for $1\le k<n$ as
the maximal value of $j$ such that $R_j\circ R_{j_{k+1}}\circ\cdots\circ
R_{j_n}(f)\ne 0$. Then
\[\gamma(f)=j_1+j_2p+j_3p^2+\cdots+j_np^{n-1}.\]
\end{proposition}

\begin{proof}
For any monomial $b_{i_1}\otimes b_{i_2}\otimes \cdots\otimes b_{i_m}$
and any $0\le j\le p-1$, we have
\[\langle b_j'|b_{i_m}\rangle\cdot b_{i_1}\otimes b_{i_2}\otimes \cdots\otimes b_{i_{m-1}}
=\left\{\begin{array}{ll}0 & \text{if $j\ne i_m$,}\\
b_{i_1}\otimes b_{i_2}\otimes\cdots\otimes b_{i_{m-1}} &
\text{otherwise.}
\end{array}\right.\]
The proof of the proposition is now straightforward.
\end{proof}

In~\cite{tullioleonovscarabottitolli}
a different basis of the dual space was used (formal dot
products with $b_k$), but the transition matrix from their basis to
$\{b_i'\}$ is triangular, so a statement similar to
Proposition~\ref{pr:heightp} holds.

One can find the height of a function $f\in\fieldp^{X^n}$ ``from the
other end'' by applying the maps $T_{b_k'}$, defined
in~\ref{ss:creationannihilation}. Recall, that these maps act by the
rule
\[T_{b_k'}(f)(x_1, x_2, \ldots, x_{n-1})=\langle b_k'|f(x, x_1, x_2,
\ldots, x_{n-1})\rangle,\]
i.e., they are linear extensions of the maps
\begin{equation}
\label{eq:Tbki1}
T_{b_k'}(b_{i_1}\otimes b_{i_2}\otimes\cdots\otimes b_{i_n})=\langle
b_k'|b_{i_1}\rangle b_{i_2}\otimes\cdots\otimes b_{i_n},
\end{equation}
and are computed by the rule
\[T_{b_k'}(f)(x_1, x_2, \ldots,
x_{n-1})=\sum_{x=0}^{p-1}\binom{x}{p-1-k}f(x, x_1, x_2, \ldots, x_{n-1}).\]

Equation~\eqref{eq:Tbki1} imply then the following algorithm for
computing height of a function.

\begin{proposition}
\label{pr:heightT}
Let $f\in\fieldp^{X^n}$, and let $h_k=T_{b_k'}(f)$ for $0\le k\le p-1$. Let $h_{k_1},
h_{k_2}, \ldots, h_{k_s}$ be the functions of the maximal height among
the functions $h_k$, and let $j_1=\max(k_i)$. Then
\[\gamma(f)=j_1+d\cdot \gamma(h_{j_1}).\]
\end{proposition}

Proposition~\ref{pr:heightT} seems to be less efficient than
Proposition~\ref{pr:heightp} in general, but it is convenient in the case
$p=2$. Let $f\in\mathbb{F}_2^{X^n}$. Denote $f_0=T_{\delta_0'}(f)$,
$f_1=T_{\delta_1'}(f)$, i.e.,
\begin{eqnarray*}f_0(x_1, x_2, \ldots, x_{n-1}) &=& f(0, x_1, x_2, \ldots, x_{n-1}),\\
f_1(x_1, x_2, \ldots, x_{n-1}) &=& f(1, x_1, x_2, \ldots, x_{n-1}).
\end{eqnarray*}

\begin{proposition}
For every $f\in\mathbb{F}_2^{X^n}$ we have
\[\gamma(f)=\left\{\begin{array}{ll}2\max(\gamma(f_0), \gamma(f_1))+1 & \text{if
      $\gamma(f_0)\ne \gamma(f_1)$,}\\ 2\gamma(f_0) & \text{if $\gamma(f_0)=\gamma(f_1)$.}\end{array}\right.\]
\end{proposition}

\begin{proof}
We have $h_0=f_1$, $h_1=f_0+f_1$.

If $\gamma(f_0)<\gamma(f_1)$, then
$\gamma(h_1)=\gamma(f_0+f_1)=\gamma(f_1)$, hence
$\gamma(f)=2\gamma(h_1)+1=2\gamma(f_1)+1$.

If $\gamma(f_0)>\gamma(f_1)$, then
$\gamma(h_1)=\gamma(f_0+f_1)=\gamma(f_0)>\gamma(f_1)=\gamma(h_0)$, hence
$\gamma(f)=2\gamma(h_1)+1=2\gamma(f_0)+1$.

If $\gamma(f_0)=\gamma(f_1)$, then $\gamma(h_1)=\gamma(f_0+f_1)<\gamma(f_1)=\gamma(f_0)$, hence $\gamma(f)=2\gamma(h_0)$.
\end{proof}

For more on height of functions on trees, and its generalizations,
see~\cite{tullioleonovscarabottitolli}.

Denote, as before, $U_n=\langle e_0, e_1, \ldots, e_n\rangle$. Then
$U_n$ consists of reduced polynomials of height not bigger than $n$.

Since the representation of $G$ is uni-triangular in the basis
$\mathsf{E}_\infty$, the spaces $U_n$ are $G$-invariant, i.e., are
sub-modules of the $G$-module $C(X^\omega, \fieldp)$. Note
also that $U_{n-1}$ has co-dimension 1 in $U_n$.

\begin{proposition}
Let $g\in\sylow$ be the adding machine. Then $U_n=(g-1)U_{n+1}$ for
every $n$.
\end{proposition}

\begin{proof}
We know that the matrix of $g$ in the basis $\mathsf{B}_\infty$ is the
infinite Jordan cell. Consequently, $(g-1)(b_0)=0$, and
$(g-1)(b_{n+1})=b_n$ for all $n\ge 0$. It follows that
\[(g-1)(U_{n+1})=(g-1)(\langle b_0, b_1, \ldots,
b_{n+1}\rangle)=\langle b_0, b_1, \ldots, b_n\rangle=U_n.\]
\end{proof}

\begin{theorem}
\label{th:uniseriality}
If $V$ is a sub-module of the $\sylow$-module $C(X^\omega, \fieldp)$,
then either $V=\{0\}$, or $V=C(X^\omega, \fieldp)$, or $V=U_n$ for
some $n$.
\end{theorem}

\begin{proof}
Let $v\in V$ and $n\ge 0$ be such that $v\in U_n\setminus
U_{n-1}$. Let $g\in\sylow$ be the adding machine defined as the
automorphism of the tree $X^*$ acting by the rule
\[g(x_1, x_2, \ldots, x_n)=\left\{\begin{array}{ll} (x_1+1, x_2,
    \ldots, x_n) & 0\le x_1\le p-2,\\
(0, g(x_2, \ldots, x_n)) & x_1=p-1.\end{array}\right.\]

Then $(g-1)^k(v)\in
U_{n-k}\setminus U_{n-k-1}$ for all $1\le k\le n$. (We assume that
$U_{-1}=\{0\}$.) It follows that $\langle v, (g-1)(v), (g-1)^2(v),
\ldots, (g-1)^n(v)\rangle_{\fieldp}=U_n\subset V$.

Let $n$ be the maximal height of an element of $V$. If $n$ is finite,
then by the proven above, $V=U_n$. If $n$ is infinite, then, by the
proven above, $V$ contains $\bigcup_{n=0}^\infty U_n=C(X^\omega, \fieldp)$.
\end{proof}

We adopt therefore, the following definition.

\begin{defi}
Let $G\le\sylow$. We say that the action of $G$ on $C(X^\omega,
\fieldp)$ is \emph{uniserial} if for
every $n\ge 0$ the set $\bigcup_{g\in G}(g-1)U_{n+1}$ generates $U_n$.
\end{defi}

A module $M$ is said to be uniserial if its lattice of sub-modules is
a chain. It is easy to see that the same arguments as in the proof
of Theorem~\ref{th:uniseriality} show that if the action of $G$ on
$C(X^\omega, \fieldp)$ is uniserial, then
$U_n$ are the only proper sub-modules of the $G$-module $C(X^\omega,
\fieldp)$. Consequently, the $G$-module $C(X^\omega, \fieldp)$ is uniserial.

In group theory (see~\cite{ddms,leedhamgreenmckay}) an action of
a group $G$ on a finite
$p$-group $U$ is said to be uni-serial, if $|N:[N, G]|=p$ for every
non-trivial $G$-invariant subgroup $N\le U$. Here $[N, G]$ is the
subgroup of $N$ generated by the elements $h^gh^{-1}$ for $h\in H$ and
$g\in G$, where $h^g$ denotes the action of $g\in G$ on $h\in H$.

Let $g\in\sylow$, and let
\[[f_0, f_1(x_1), f_2(x_1, x_2), \ldots]\]
be the tableau of $g$. We have seen in Subsection~\ref{ss:principal}
(see the proof of Proposition~\ref{pr:principal})
that the entries of the principal columns $(a_{0, p^n}, a_{1, p^n},
\ldots, a_{p^n-1, p^n})^\top$ of the matrix $(a_{i, j})_{i,
  j=0}^\infty$ of $\pi_\infty(g)$ in the basis $\mathsf{E}_\infty$ are
precisely the coefficients of the polynomials $f_n$:
\[f_n(x_1, x_2, \ldots, x_{n-1})=\sum_{k=0}^{p^n-1}a_{k, p^n}e_k,\]
where $e_k$ is the monomial of height $k$.

It follows that the height of $f_n$ is equal to the largest index of
a non-zero non-diagonal entry of the column number $p^n$ of the matrix
of $\pi_\infty(g)$ in the basis $\mathsf{E}_\infty$.
Note that the same is true for the matrix of
$\pi_\infty(g)$ in the basis $\mathsf{B}_\infty$.

\begin{proposition}
\label{pr:uniserialitycriterion}
Let $G\le\sylow$, and let
$\alpha:\sylow\arr\fieldp^\omega:g\mapsto(\alpha_0(g), \alpha_1(g), \ldots)$ be the
abelianization homomorphism given by~\eqref{eq:alphamap}. The action of the group $G$ on
$C(X^\omega, \fieldp)$ is uniserial if and only
if every homomorphism $\alpha_k:G\arr\fieldp$ is non-zero.
\end{proposition}

\begin{proof}
It follows from Theorem~\ref{th:diagonal} that all homomorphisms
$\alpha_k$ are non-zero if and only if for every $k=1, 2, \ldots$
there exists $g_k\in G$ such that the entry number $k$ on the first
diagonal of $\pi_\infty(g_k)$ is non-zero.

Then for every monomial $e_k$ the height of $(1-g_k)(e_k)$ is equal to
$k-1$, which shows that $\bigcup_{i=1}^k(1-g_k)(U_k)$ generates
$U_{k-1}$, hence the action of $G$ is uniserial.
\end{proof}

\begin{corollary}
Let $S$ be a generating set of $G\le\sylow$. Then the action of $G$ on
$C(X^\omega, \fieldp)$ is uniserial if
and only if for every $k=0, 1, \ldots$ there exists $g_k\in S$ such
that $\alpha_k(g_k)\ne 0$.
\end{corollary}

Note that it also follows from Theorem~\ref{th:diagonal} and from the
fact that the entries in the principal columns are the coefficients of
the polynomials in the tableau, that
$\alpha_n(g)\ne 0$ if and only if height of the polynomial $f_n$ of
the tableau $[f_0, f_1(x_1), f_2(x_1, x_2), \ldots]$ representing $g$
is equal to $p^n-1$, i.e., has the maximal possible value.

\begin{examp}
The cyclic group generated by an element $g\in\sylow$ is transitive on
the levels $X^n$ if and only if $\alpha_n(g)\ne 0$ for all $n$. It
follows that if $G$ contains a level-transitive element, then its
action is uniserial. But there exist torsion groups with uniserial
action on $C(X^\omega, \fieldp)$, as the following
example shows.
\end{examp}

\begin{examp}
It is easy to check that for the generators $a, b, c, d$ of the
Grigorchuk group, we have $\alpha(a)=(1, 0, 0, 0, \ldots)$, and
\begin{eqnarray*}
\alpha(b) &=& (0, 1, 1, 0, 1, 1, 0, 1, 1, 0, \ldots)\\
\alpha(c) &=& (0, 1, 0, 1, 1, 0, 1, 1, 0, 1, \ldots)\\
\alpha(d) &=& (0, 0, 1, 1, 0, 1, 1, 0, 1, 1, \ldots).
\end{eqnarray*}
(In the last three equalities, each sequence have a pre-period of length
1 and a period of length 3.) It follows that the action of the Grigorchuk group is uniserial.
\end{examp}

\begin{examp}
Gupta-Sidki group~\cite{gupta-sidkigroup} is generated by two elements $a, b$
acting on $\{0, 1, 2\}^*$, where $a$ is the cyclic permutation
$\sigma=(012)$ on the first level of the tree (i.e., changing only the
first letter of a word), and $b$ is defined by the wreath recursion
\[b=(a, a^{-1}, b).\]
Then $\alpha(a)=(1, 0, 0, \ldots)$, and $\alpha(b)=(0, 0, 0, \ldots)$,
hence the group $\langle a, b\rangle$ does not act uniserially on
$\{0, 1, 2\}^\omega$.
\end{examp}

\def\cprime{$'$}
\providecommand{\bysame}{\leavevmode\hbox to3em{\hrulefill}\thinspace}
\providecommand{\MR}{\relax\ifhmode\unskip\space\fi MR }
\providecommand{\MRhref}[2]{%
  \href{http://www.ams.org/mathscinet-getitem?mr=#1}{#2}
}
\providecommand{\href}[2]{#2}

\end{document}